\documentclass[amssymb,twoside,12 pt,amscd,reqno]{amsart}
\usepackage{amsmath}
\usepackage{amscd}
\usepackage{latexsym}
\theoremstyle{plain}
\newtheorem{thm}{Theorem}[section]
\newtheorem{cor}[thm]{Corollary}
\newtheorem{lem}[thm]{Lemma}
\newtheorem{prop}[thm]{Proposition}
\newtheorem{Claim}[thm]{Claim}

\theoremstyle{definition}
\newtheorem{defn}[thm]{Definition}

\newcommand{\Hdg}[2]{\mathcal{H}^{#1,#2}}

\newcommand{\ZZ}{\mathbb{Z}}
\newcommand{\CC}{\mathbb{C}}

\newcommand{\AAA}{\mathbb{A}}
\newcommand{\PP}{\mathbb{P}}

\newcommand{\HI}[2]{\text{Hilb}_{#1}(#2)}

\newcommand{\GR}[3]{\text{Grass}_{#1}(#2,#3)}

\def\ni{\noindent}

\def\pic{{\rm Pic}}

\def\P{{\mathbb P}}

\def\Z{{\mathbb Z}}

\def\ps{\vspace{4pt}}

\newcommand{\lt}{\left}
\newcommand{\rt}{\right}

\newcommand{\OO}{\mathcal O}

\begin{document}

\title{Abel-Jacobi Maps Associated to Smooth Cubic Threefolds}
\author[Harris]{Joe Harris} 
\address{Department of Mathematics \\ Harvard University \\ Cambridge MA 02138}
\email{harris@math.harvard.edu}
\author[Roth]{Mike Roth} 
\address{Department of Mathematics \\ University of Michigan \\ Ann
Arbor, MI 48109}  
\email{mikeroth@math.lsa.umich.edu}
\author[Starr]{Jason Starr}\thanks{The third author was partially
supported by an NSF Graduate Research Fellowship and a Sloan
Dissertation Fellowship} 
\address{Department of Mathematics \\ Massachusetts Institute of Technology \\ Cambridge MA 02139}
\email{jstarr@math.mit.edu}
\date{\today}

\maketitle

\tableofcontents

\begin{abstract}
In this article we consider the spaces $\Hdg{d}{g}(X)$ parametrizing
curves of degree $d$ and genus $g$ on a smooth cubic threefold $X
\subset \P^4$, with regard in particular to the Abel-Jacobi map $u_d :
\Hdg{d}{g}(X) \to J^3(X)$ to the intermediate Jacobian $J^3(X)$ of
$X$.  Our principle result is that for all $d \leq 5$ the map $u_d$
coincides with the maximal rationally connected fibration of
$\Hdg{d}{g}(X)$.
\end{abstract}

\section{Introduction} \label{sec-intro}

In this paper we will study spaces parametrizing curves on a smooth
cubic threefold $X \subset \P^4$. There has been a great deal of work
in recent years on the geometry of spaces parametrizing rational
curves on a variety $X$.  Most of it has focussed on the enumerative
geometry of these spaces: the description of their Chow rings and the
evaluation of certain products in their Chow rings. Here we will be
concerned with a very different sort of question: we will be concerned
with the birational geometry of the spaces $M$.

\ps

We should start by explaining, at least in part, our motivation. The
central object in curve theory---the one that links together every
aspect of the theory, and whose study yields the majority of theorems
in the subject---is the Abel-Jacobi map. This is the map from the
symmetric product $C_d$ of a curve $C$, parametrizing $0$-cycles of a
given degree $d$ on $C$, to the Jacobian $\pic^d(C) \cong J(C)$,
defined variously as the space of cycles of degree $d$ mod linear
equivalence (that is, the space of line bundles of degree $d$ on $C$)
or, over the complex numbers, as the complex torus
$H^1(C,K)^*/H_1(C,\Z)$.

\ps

Under the circumstances, it's natural to ask what sort of analogue of
this we might be able to define and study in higher dimensions. There
have been numerous constructions proposed to this end; they have been
too varied to categorize easily, but the majority adopt one of two
approaches, corresponding to the definitions of $\pic^d(C) \cong J(C)$
for a curve. In the first, we look again at the space of cycles on a
variety $X$ mod some equivalence relation (typically rational
equivalence); in the second, we try to form a geometric object out of
the Hodge structure of $X$.

\ps

Both approaches suffer from some seemingly unavoidable
difficulties. In the first, the quotients of the spaces on cycles on a
variety $X$ by rational equivalence tends to be either too big or too
small: if we mod out by algebraic equivalence we tend to lose too much
information, while if we mod out by rational equivalence the quotient
is too large (and in particular too hard to calculate even in simple
concrete cases).  As for the second, we simply don't know how to make
an algebrao-geometric object out of a Hodge structure except in very
special cases---specifically, the case of Hodge structures of odd
weight with all but two of the Hodge groups vanishing, where we can
form the intermediate Jacobian and it will be an Abelian variety. In
these cases the theory has had spectacular successes (such as the
proof of the irrationality of cubic threefolds), but it's undeniably
frustrating that it has no wider application.

\ps

In this paper we'd like to suggest another approach to the problem,
based on a construction of Koll\'ar, Miyaoka and Mori: the {\it
maximal rationally connected fibration} associated to a variety
$X$. Briefly, we say that a smooth projective variety $X$ is {\it
rationally connected} if two general points of $X$ can be connected by
a chain of rational curves. This condition is equivalent to
rationality and unirationality for curves and surfaces; in higher
dimensions it is weaker than either.  More generally, the maximal
rationally connected fibration associates to a variety $X$ a
(birational isomorphism class of) variety $Z$ and a rational map $\phi
: X \to Z$ with the properties that
\begin{itemize}
\item the fibers $X_z$ of $\phi$ are rationally connected; and conversely 
\item almost all the rational curves in $X$ lie in fibers of $\phi$:
for a very general point $z \in Z$ any rational curve in $X$ meeting
$X_z$ lies in $X_z$.
\end{itemize}
\ni The variety $Z$ and morphism $\phi$ are unique up to birational 
isomorphism, and are called the
\textit{mrc quotient} and \textit{mrc fibration} of $X$, respectively.
They measure the failure of $X$ to be rationally connected: if $X$ is
rationally connected, $Z$ is a point, while if $X$ is not uniruled we
have $Z=X$.

\ps

The point of this is, if we start with the symmetric product $C_d$ of
a curve $C$ and apply this construction, we find (for $d$ large, at
any rate) that the mrc fibration of $C_d$ is the Abel-Jacobi map $\phi
: C_d
\to J(C)$. In other words, we don't have to define the Jacobian 
$J(C)$ either in terms of cycles mod linear
equivalence or via Hodge theory; we can realize it simply the mrc
quotient of the space $C_d$ parametrizing $0$-cycles of reasonably
large degree $d$ on $C$. The obvious question is then: what happens
when we take the mrc quotient of the variety parametrizing
cycles\footnote{Note that we have a choice of parameter spaces for the
cycles on $X$: the Chow variety, the Hilbert scheme and, in the case
of 1-dimensional cycles, the Kontsevich space. But since we are
concerned with the birational geometry of these spaces it really
doesn't matter which we choose to work with, except for technical
convenience.} on a higher-dimensional variety?

\ps

To begin with, taking the mrc quotient of the varieties parametrizing
$0$-cycles on $X$ doesn't seem to yield much: if $X$ is rationally
connected, so will be its symmetric powers; and if $X$ possesses any
holomorphic forms one expects the dimensions of the mrc quotients of
$X_d$ tend  
to $\infty$ with $d$ (c.f. ~\cite{RE}).  We
turn our attention next to curves on $X$ and ask: what can we say
about the mrc quotients of the varieties $\Hdg{d}{g}(X)$ parametrizing
curves of degree $d$ and genus $g$ on $X$? For various reasons it
makes sense to look primarily at the spaces parametrizing rational
curves---for one thing, we don't want to get involved in the geometry
of the moduli space $M_g$ of curves of genus $g$ when we are only
interested in invariants of $X$---but in the course of studying
rational curves we will also discover facts about the geometry of
$\Hdg{d}{g}(X)$ for $g > 0$.

\ps

For $X = \P^n$ and for $X \subset \P^{n+1}$ a quadric hypersurface,
the answer is known, at least in the case of rational curves, and is
trivial: the variety $\Hdg{d}{0}(X)$ parametrizing rational curves of
degree $d$ on $X$ is itself rationally connected in both cases. The
first real test is thus cubic hypersurfaces, and here we will focus
specifically on the case of a smooth cubic threefold $X
\subset
\P^4$.

\ps

We start by asking the most naive possible question: is the variety
$\Hdg{d}{g}(X)$ rationally connected? The answer is that it can never
be: we have the Abel Jacobi map $$ u : \Hdg{d}{g}(X) \to J^3(X) $$
from the Kontsevich space to the intermediate Jacobian $$ J_3(X) =
H^{2,1}(X)^*/H_3(X,\Z).  $$ Since any map of a rational curve to a
complex torus/Abelian variety is constant, any rational curve on
$\Hdg{d}{g}(X)$ must lie in a fiber of $u$. The Abel-Jacobi map thus
factors through the mrc fibration $\phi : \Hdg{d}{g}(X) \to Z$ of
$\Hdg{d}{g}(X)$.

\ps

Given this, the second most naive question would be: is the
intermediate Jacobian $J^3(X)$ the mrc quotient of $\Hdg{d}{g}(X)$?
This may seem equally naive, but it turns out to hold for small values
of $d$. The main result of this paper is the

\ps

\begin{thm}\label{mainthm}  If $X \subset \P^4$ is any smooth 
cubic threefold, for any $d \leq 5$ the Abel-Jacobi
map $$ u : \Hdg{d}{g}(X) \to J^3(X) $$ is the maximal rationally
connected fibration of $\Hdg{d}{g}(X)$.
\end{thm}

As we indicated, we view this result as simply the first test of a
general hypothesis. It would be very interesting to see if the
analogous statement holds for other rationally connected threefolds.
Work is being done in this direction: Ana-Maria Castravet has proved
in ~\cite{C} that for $X \subset \P^5$ the intersection of two quadrics, the
analogous statement---that the mrc quotient of the space
$\Hdg{d}{0}(X)$ is the intermediate Jacobian---holds for all $d$. And
work of Harris and Starr~\cite{HS} shows that the same is true at least
for components of the variety $\Hdg{d}{0}(X)$ corresponding to an open
cone of curve classes $\beta \in H_2(Z,\Z)$ in a blow-up $X$ of
$\P^3$.

\ps

Ultimately, though, we want to extend this research to
higher-dimensional varieties $X$, where the Hodge theory doesn't seem
to provide us with a candidate for the mrc quotients of the spaces of
rational curves on $X$. To relate this to a specific issue: for smooth
cubic threefolds $X$, the fact that the Hodge structure on $H^3(X,\Z)$
is not isomorphic to a product of Hodge structures of curves implies
that $X$ is not rational. A close examination of the Hodge structure
of a very general cubic fourfold $X \subset \P^5$ suggests a similar
finding: the Hodge structure on $H^4(X,\Z)$ does not appear to be a
product of factors of Hodge structures of surfaces (see [Has]), and if
this is indeed the case it would imply that $X$ is irrational. But the
fact that the Hodge structure is not conveyed in the form of a
geometric object---there is no construction analogous to that of the
intermediate Jacobian of a cubic threefold---has frustrated our
attempts to make this into a proof. But now we can ask: if the
intermediate Jacobian of a cubic threefold $X$ may be realized as the
mrc quotient of the space of rational curves on $X$, what happens when
we take the mrc quotient of the space of rational curves on a cubic
fourfold?

\ps

In this paper, though, we will be concerned exclusively with the
geometry of cubic threefolds. We start some preliminary sections on
the Abel-Jacobi map for curves on threefolds in general and on the
geometry of curves of low degree on cubic threefolds in particular. We
then launch into the analysis of the Abel-Jacobi map for curves of
degree 3, 4 and 5 on a smooth cubic threefold $X$. The basic technique
here is {\it residuation}: we associate to a curve $C \subset X$ on
$X$ and a surface $S \subset \P^4$ containing it, a residual curve
$C'$ with (roughly) $C +C' = S
\cap X$. 
In this way we relate the parameter space for curves $C$ to that for
curves $C'$.
In ~\cite{HRS1} it is proved
that for $d\leq 5$, each of the spaces $\Hdg{d}{g}(X)$ is
irreducible of dimension $2d$.  Combining this with the residuation
technique, we prove that for $d\leq 5$ the fibers of the Abel-Jacobi
map $\Hdg{d}{g}(X)\to J^3(X)$ are unirational.

\ps

Many of the results in this paper have also been proved in ~\cite{MT}
and ~\cite{AM} by considering moduli of vector bundles on $X$.  Our
theorem~\ref{thm-st30} corresponds to \cite[proposition 3.2]{MT}.  Our
corollary~\ref{cor-51fib} corresponds to ~\cite[theorem 5.6(iii)]{MT}
and ~\cite[theorem 3.2]{AM}.  Our proof is original and does not use
vector bundles on $X$, but rather an analysis of the divisors in the
intermediate Jacobian.  In addition, we give an analysis of quartic
rational curves (which are also considered in ~\cite{MT}), of quintic
rational curves, and of quintic curves of genus $2$.  Furthermore,
we prove that for every degree $d>3$, the
Abel-Jacobi map $\alpha_{d,0}:\Hdg{d}{0}(X)\rightarrow
J(X)$ is dominant and the general fiber is irreducible.

\subsection{Notation} \label{sec-not}

All schemes in this paper will be schemes over $\CC$.  All absolute
products will be understood to be fiber products over
$\text{Spec}(\CC)$.

\ps

For a projective variety $X$ and a numerical polynomial $P(t)$,
$\HI{P(t)}{X}$ denotes the corresponding Hilbert scheme.  For integers
$d,g$, $\Hdg{d}{g}(X)\subset \HI{dt+1-g}{X}$ denotes the open
subscheme parametrizing smooth, connected curves of degree $d$ and
genus $g$, and $\overline{\Hdg{d}{g}}(X)$ denotes the closure of
$\Hdg{d}{g}(X)$ in $\HI{dt+1-g}{X}$.

\section{Review of the Abel-Jacobi Map} \label{sec-rev}

Our object of study are the \emph{Abel-Jacobi maps} associated 
to families of 1-cycles on a smooth cubic hypersurface $X\subset
\PP^4$.  The reader is referred to ~\cite{CG} and ~\cite{VHS}
for full definitions.  Here we recall only a few facts
about Abel-Jacobi maps.

\ps

Associated to a smooth, projective threefold $X$ there is a complex
torus  
\begin{equation}
J^2(X) = H^3_\ZZ(X)\backslash H^3(X,\CC) / \lt(H^{3,0}(X)\oplus
H^{2,1}(X)\rt).
\end{equation}
In case $X$ is a cubic hypersurface in $\PP^4$ (in fact for any
\emph{rationally connected} threefold) then $J^2(X)$ is a principally
polarized
abelian variety with theta divisor $\Theta$.  Given an algebraic
1-cycle $\gamma\in A_1(X)$ which 
is \emph{homologically equivalent to zero}~\cite[13]{VHS}, one can
associate a point $u_2(\alpha)$.  The construction is analogous to
the Abel-Jacobi map for a smooth, projective algebraic curve $C$ which 
associates to each 0-cycle $\gamma\in A_0(C)$ which is homologically
equivalent to zero a point $u_1(\alpha)\in J^1(C)$, the Jacobian
variety of $C$.  In particular $u_2:A_1(X)^{hom}\rightarrow J^2(X)$ is 
a group homomorphism.

\ps

Suppose that $B$ is a normal, connected variety of dimension $n$ and
$\Gamma\in A_{n+1}(B\times X)$ is an $(n+1)$-cycle such that for each
closed point $b\in B$ the corresponding cycle $\Gamma_b\in
A_1(X)$~\cite[\S 10.1]{F} is homologically equivalent to zero.  Then
in this case the set map $b\mapsto u_1(\Gamma_b)\in J^2(X)$ comes from 
a (unique) algebraic morphism $u=u_\Gamma:B\rightarrow J^2(X)$.  We
call this morphism the \emph{Abel-Jacobi map} determined by $\Gamma$.

\ps

More generally, suppose $B$ as above, $\Gamma\in A_{n+1}(B\times
X)$ is any $(n+1)$-cycle, and suppose $b_0\in B$ is some base-point.
Then we can form a new cycle $\Gamma'=\Gamma - \pi_2^*\Gamma_{b_0}$,
and for all $b\in B$ we have $\Gamma'_b = \Gamma_b - \Gamma_{b_0}$ is
homologically equivalent to zero.  Thus we have an algebraic morphism
$u=u_{\Gamma'}:B\rightarrow J^2(X)$.  Of course this morphism depends
on the choice of a base-point, but changing the base-point only
changes the morphism by a constant translation.  Thus we shall speak
of any of the morphisms $u_{\Gamma'}$ determined by $\Gamma$ and the
choice of a base-point as an \emph{Abel-Jacobi map} determined by
$\Gamma$.  

\ps

Suppose that $\Gamma_1,\Gamma_2\in A_{n+1}(B\times X)$ are two
$(n+1)$-cycles.  Then $u_{\Gamma_1+\Gamma_2}$ is the pointwise sum
$u_{\Gamma_1}+u_{\Gamma_2}$.  This trivial observation is frequently
useful.  

\ps 

\textbf{The Residuation Trick:}
Another useful observation is that any Abel-Jacobi morphism
$\alpha_\Gamma$ contracts all rational curves on $B$, since an Abelian 
variety contains no rational curves.  Combined with the observation in
the last paragraph, this leads to the residuation trick:  Suppose that
$B$ is a normal, rationally-connected variety and $\Gamma\in A_{n+1}(B\times 
X)$ is an $(n+1)$-cycle.  Then $u_\Gamma:B\rightarrow J(X)$ is a
constant map.  Now suppose that $B'\subset B$ is a normal closed
subvariety and that $\Gamma|_{B'\times X}$ decomposes as a sum of
cycles $\Gamma_1+\Gamma_2$.  Since $u_{\Gamma_1}+u_{\Gamma_2}$ equals
a constant map, we conclude that $u_{\Gamma_1}$ is the pointwise
inverse of $u_{\Gamma_2}$, up to a fixed additive constant.

\section{Lines, Conics and Plane Cubics}~\label{sec-123}

The study of the Abel-Jacobi map associated to the space
$\Hdg{1}{0}(X)$ of lines on $X$ was carried out in  ~\cite{CG}.  
In this section we will summarize their results, which will also be
useful to us for studying Abel-Jacobi maps of higher degree curves.
In this section we also consider the Abel-Jacobi maps associated to
the spaces $\Hdg{2}{0}(X)$ and $\Hdg{3}{1}(X)$ of plane conics and
plane cubics in $X$.  In each case the Abel-Jacobi is trivial to
describe.  

\subsection{Lines}\label{subsec-lines}

For brevity we refer to the \emph{Fano scheme of lines},
$\Hdg{1}{0}(X)$, simply as $F$.  
Two general lines $L_1,L_2\subset \PP^4$ determine a hyperplane by
$\text{span}(L_1,L_2)$.  We generalize this as follows:
Let $(F\times F-\Delta)\xrightarrow{\Phi} 
\PP^{4\vee}$ denote the following set map:
\begin{equation}\label{eqn-16}
 \Phi \left( [L_1,L_2] \right) =
  \left\{ \begin{array}{ll}
    \lt[ \text{span} \left( L_1,L_2 \right)\rt] &\text{if }L_1\cap L_2
  = \emptyset,  \\
  \lt[T_p X\rt]  & \text{if }p\in L_1\cap L_2
 \end{array} \right.
\end{equation} 
By~\cite[lemma 12.16]{CG}, $\Phi$ is algebraic.  
Let $X^\vee\subset \PP^{4\vee}$ denote the dual variety of $X$,
i.e. the variety 
parametrizing tangent hyperplanes to $X$.  Let $X^\vee_s\subset
X^\vee$ denote 
the subvariety parametrizing hyperplanes $H$ which are tangent to $X$ and
such that the singular locus of $H\cap X$ is not simply a single ordinary 
double point.  Let $U_s\subset U\subset F\times F$ denote the open
sets $\Phi^{-1}\lt(\PP^{4\vee} - X^\vee\rt)\subset \Phi^{-1}\lt(\PP^{4\vee} - 
X^\vee_s\rt)$.  And let $I\subset F\times F$ denote the divisor
parametrizing 
incident lines, i.e. $I$ is the closure of the set $\{\lt([L_1],[L_2]\rt) 
: L_1\neq L_2, L_1\cap L_2\neq\emptyset\}$.
In~\cite{CG}, Clemens and Griffiths completely describe both
the total 
Abel-Jacobi map $F\times F\xrightarrow{\psi} J(X)$ and the
Abel-Jacobi map $F\xrightarrow{i} J(X)$.  Here is a summary of their results

\ps

\begin{thm} \label{thm-cg}
\begin{enumerate} 

\item The Fano variety $F$ is a smooth surface and the 
Abel-Jacobi
map $F\xrightarrow{u} J(X)$ is a closed immersion~\cite[theorem 7.8,
theorem 12.37]{CG}. 

\item The induced map $Alb(F)=J^2(F)\rightarrow J(X)$ is an isomorphism of 
principally polarized Abelian varieties~\cite[theorem 11.19]{CG}.
\item The class of $u(F)$ in $J(X)$
is $\frac{[\Theta]^3}{3!}$~\cite[proposition 13.1]{CG}.

\item The difference of 
Abel-Jacobi maps 
\begin{equation}
\psi:F\times F \rightarrow J(X), \ \ \psi([L],[L']) = u([L]) - u([L'])
\end{equation}
maps $F\times F$ generically $6$-to-$1$ to the theta divisor
$\Theta~\subset ~J(X)$~\cite[section 13]{CG}.

\item Let $(\Theta-\{0\}) \xrightarrow{\mathcal{G}} 
\PP H^{1,2}(X)^\vee$ denote the Gauss map.  If we identify
$\PP H^{1,2}(X)$ with $\PP^4$ via the
Griffiths residue calculus~\cite{GR}, then the composite map
\begin{equation}
(F\times_C F-\Delta)\xrightarrow{\psi}(\Theta -\{0\})\xrightarrow{\mathcal{G}}
\PP^{4\vee}
\end{equation}
is just the map $\Phi$ defined above~\cite[formula 13.6]{CG}.

\item The fibers of the Abel-Jacobi map form a \emph{ Schl\"afli
double-six}, i.e. the general fiber of $\psi:F\times F\rightarrow J$
is of the form $\{(E_1,G_1),\dots,(E_6,G_6)\}$ where the lines
$E_i,G_j$ lie in a smooth hyperplane section of $X$, the $E_i$ are
pairwise skew, the $G_j$ are pairwise skew, and $E_i$ and $G_j$ are
skew iff $i=j$.  

There is a more precise result than above.  Let 
\begin{equation}R' \subset \lt(U\times_{\PP^{4\vee}} U\rt) \times F \times
\GR{}{3}{V}\times \GR{}{3}{V}\end{equation} 
be the closed subscheme parametrizing data
$\lt(\lt([L_1],[L_2]\rt),\lt([L_3],[L_4]\rt),[l], [H_1], [H_3]\rt)$
such that for each $i=1,\dots,4$, $l\cap L_i\neq\emptyset$ and such
that $H_1\cap X = l\cup L_1\cup L_4$, $H_2\cap X = l\cup L_2\cup L_3$.
Let $R\subset U\times_{\PP^{4\vee}} U$ be the image of $R'$ under the
projection map.  Let $\Delta\subset U\times U$ be the diagonal.  Then
the fiber product $U\times_\Theta U\subset U\times U$ is just the
union $R\cup\Delta$~\cite[p. 347-348]{CG}

\item The branch locus of $\Theta\xrightarrow{\mathcal{G}} \PP^{4\vee}$ equals
the branch locus of $F\times F \xrightarrow{\Phi} \PP^{4\vee}$ equals   
the dual variety of $X$, i.e. the variety 
parametrizing the tangent hyperplanes to $X$.  The 
ramification locus of $U \xrightarrow{\Phi} \PP^{4\vee}$ equals the
ramification 
locus of $U\xrightarrow{\psi} \Theta$ equals
the divisor $I$.  Each
such pair is a simple ramification point of both $\psi$ and
$\Phi$~\cite[lemma 13,8]{CG}.  
\end{enumerate}
\end{thm}

In~\cite{Tj}, Tjurin also analyzed cubic threefolds and the associated 
intermediate Jacobians.  We summarize his results:

\ps

\begin{thm}~\cite{Tj} \label{thm-tj}  
Let $\widetilde{J(X)}$ be the variety obtained by blowing up $0\in J(X)$ and
let $\widetilde{\Theta}$ be the proper transform of $\Theta$.  Let 
$\widetilde{F\times F}$
be the variety obtained by blowing up the diagonal in $F\times F$.  The
exceptional divisor $E\subset\widetilde{J(X)}$ is isomorphic to $\PP^4$ and
this isomorphism identifies the intersection $\widetilde{\Theta}\cap E$ with
our original cubic threefold $X$.  The exceptional divisor $E'\subset
\widetilde{F\times F}$ is the projective bundle $\PP T_F$.  
In fact there is
an isomorphism of sheaves of the tautological rank 2 sheaf $S=S(2,V)$
and the tangent bundle $T_F$ so that $E'\cong \PP_{F}S$.  The induced
morphism $\PP_{F}S\xrightarrow{\tilde{\psi}} \PP^4$ is just the usual
morphism induced by the map of sheaves $S\hookrightarrow V\otimes_\CC \OO_F$.
In particular, this morphism is generically $6$-to-$1$ onto $X$ with 
ramification locus $\Sigma_1\subset F$ consisting of the lines of "second
type", i.e. lines $L\subset X$ such that $N_{L/X}\cong \OO_L(-1)\oplus\OO_L(1)$.
Let $D\subset F\times F$ denote the divisor of intersecting lines and
let $\widetilde{D}$ be the proper transform of $D$.  Then $\widetilde{D}\cap E'
=\Sigma_1$ so that we finally have the result: The morphism
$\widetilde{F\times F}\rightarrow{\widetilde{\Theta}}$ is finite of degree
$6$ with ramification locus $\widetilde{D}$.
\end{thm}

\subsection{Conics}\label{subsec-cons}

Next we consider the space $\Hdg{2}{0}(X)$ parametrizing smooth plane
conics on $X$.  Given a plane conic $C\subset X$, consider the plane
cubic curve which is the intersection $\text{span}(C)\cap X$.  This
curve contains $C$ is an irreducible component, and the residual
component is a line $L\subset X$.  In this way we have a morphism
$\Hdg{2}{0}(X) \rightarrow F$ by $[C]\mapsto [L]$.  Conversely, given
a line $L\subset X$ and a $2$-plane $P\subset \PP^4$ which contains
$L$, then the residual to $L$ in $P\cap X$ is a plane conic.  In this
way one sees that $\Hdg{2}{0}(X)$ is isomorphic to an open subset of
the $\PP^2$-bundle
\begin{equation}
\PP Q = \{([L],[P])\in F\times \mathbf{G}(2,4) | L\subset P\}.
\end{equation}
In ~\cite[prop 3.4]{HRS1} we prove the stronger result:

\begin{prop}~\label{prop-cons1}  
The morphism
$\Hdg{2}{0}(X)\rightarrow F$ is isomorphic to an open subset of a
$\PP^2$-bundle $\PP Q\rightarrow F$.  In particular, $\Hdg{2}{0}(X)$
is smooth and connected of dimension $4$.  Moreover $\PP Q$ is the
normalization of $\HI{2t+1}{X}$.  
\end{prop}

Next we describe the Abel-Jacobi map for conics.  
By the residuation trick, the Abel-Jacobi map $u_{D_2}:\PP
Q\rightarrow J(X)$ is, up to a 
fixed translation, the pointwise inverse of the composite $\PP
Q\xrightarrow{\pi} F\xrightarrow{u} J(X)$.  Thus the fibers of
$u_{D_2}:\PP Q\rightarrow J(X)$ are just the fibers of $\pi$,
i.e. $\PP^2$'s.

\subsection{Plane Cubics}\label{sec-pcubs}

Every curve $C\subset \PP^4$ with Hilbert polynomial $3t$ is a plane
cubic, and the 2-plane $P=\text{span}(C)$ is unique; we have
that $C=X\cap P$.  Therefore the Hilbert scheme $\HI{3t}{X}$ is just
the Grassmannian $\mathbb{G}(2,\PP^4)$ of 2-planes in $\PP^4$ and
$\Hdg{3}{1}(X)$ is just an open subset of $\mathbb{G}(2,\PP^4)$.    
Since
$\mathbb{G}(2,\PP^4)$ is rational, it follows that both the total Abel-Jacobi
map $\psi:\HI{3t}{X}\times \HI{3t}{X}\rightarrow J(X)$ and the
Abel-Jacobi map $i:\HI{3t}{X}\rightarrow J(X)$ are constant maps.  

\subsection{Two Useful Results}\label{sec-res}

Recall that $\psi:F\times F\rightarrow J(X)$ is defined as the
\emph{difference map} $\psi = u\circ \pi_1-u\circ \pi_2$.  Now we
consider the sum map $\psi'=u\circ\pi_1+u\circ\pi_2$.  
Alternatively one may define $\psi'$ to be the (Weil extension of the)
morphism $(F\times F)-\Delta \rightarrow J(X)$ which is the
Abel-Jacobi map corresponding to the flat family $Z\xrightarrow{\pi}
(F\times F-\Delta)$ defined as follows.  For $L_1$ and $L_2$ skew
lines, we define $Z_{\lt([L_1],[L_2]\rt)}$ to be the disjoint union of
$L_1$ and $L_2$.  If $p\in L_1\cap L_2$, we define
$Z_{\lt([L_1],[L_2]\rt)}$ to be the subscheme whose reduced scheme is
just the union $L_1\cup L_2$ and which has an embedded point at $p$
corresponding to the normal direction of $\text{span}(L_1,L_2)\subset
T_p X$.

\ps

Clearly $\psi'$ commutes with the involution 
\begin{equation}
\iota:F\times F \rightarrow F\times F, ([L_1],[L_2])\mapsto
([L_2],[L_1]). 
\end{equation}
So $\psi'$ factors through the
quotient $F\times F\rightarrow \text{Sym}^2 F$.  We will denote by
$\text{Sym}^2 F
\xrightarrow{\psi''} J(X)$ the induced morphism.
Let us define $\Theta'$ to be the scheme-theoretic image of $\psi'$.

\ps

First we prove an easy lemma about
cohomology of complex tori.  Given a $g$-dimensional complex torus $A$
and a sequence of nonzero integers $(n_1,\dots,n_r)$, consider the
holomorphic map $f=f_{(n_1,\dots,n_r)}:A^r\rightarrow A$ defined by
$f(a_1,\dots,a_r) = n_1a_1 + \dots + n_r a_r$.  There is an induced
Gysin image map on cohomology $f_*:H^{p+2(r-1)g}(A^r)\rightarrow
H^p(A)$.  Now by the K\"unneth formula, the cohomology $H^q(A^r)$ is a
direct sum of K\"unneth components $H^{q_1}(A)\otimes \dots \otimes
H^{q_r}(A)$ with $q=q_1+\dots+q_r$.

\begin{lem}\label{lem-cohom}  For each integer $r\geq 1$, for each
integer $0\leq p\leq 2g$, and for each decomposition $(q_1,\dots,q_r)$
with $q_1+\dots+q_r = p + 2(r-1)g$, there is a homomorphism
\begin{equation}
g=g_{(q_1,\dots,q_r)}:H^{q_1}(A)\otimes\dots\otimes
H^{q_r}(A)\rightarrow H^p(A)
\end{equation}
such that for each sequence of nonzero integers $(n_1,\dots,n_r)$,
with $f$ defined as above, the restriction of the Gysin image map
\begin{equation}
f_*:H^{q_1}(A)\otimes \dots\otimes H^{q_r}(A)\rightarrow H^p(A)
\end{equation}
satisfies $f_* = \lt(n_1^{2g-q_1}\cdot \dots \cdot n_r^{2g-q_r}\rt) g$.  
\end{lem}

\begin{proof}  This is just a computation.  By Poincar\'e duality, to
give the homomorphism $g$, it is equivalent to give a bilinear pairing
\begin{equation}
\lt(H^{q_1}(A)\otimes \dots \otimes H^{q_r}(A)\rt) \times H^{2g-p}(A)
\rightarrow \ZZ. 
\end{equation}
Moreover, since $H^{2g-p}(A)= \bigwedge^{2g-p}H^1(A)$, it suffices to
define the pairing for pure wedge powers
$\alpha = \alpha_1\wedge\dots\wedge\alpha_{2g-p}\in H^{g-p}(A)$.  
Define $S$ to be the set of functions
$\sigma:\{1,2,\dots,2g-p\}\rightarrow {1,2,\dots, r}$ such that for each 
$i=1,\dots, r$, we have $q_i+\#\sigma^{-1}(i) = 2g$.  
We define the pairing by taking the wedge product in $H^*(A^r)$ and
then taking the degree as follows:
\begin{equation}
\langle \beta = \beta_1\otimes \dots \otimes \beta_r, \alpha_1\wedge\dots
\wedge \alpha_{2g-p}\rangle = \text{deg}\sum_{\sigma\in S}
\pi_1^*\beta_1 \wedge \dots \wedge \pi_r^*\beta_r \wedge
\pi_{\sigma(1)}^*\alpha_1 \wedge \dots \wedge
\pi_{\sigma(2g-p)}^*\alpha_{2g-p}.
\end{equation}

\

The fact that $f_* = \lt(n_1^{2g-q_1}\cdot \dots \cdot
n_r^{2g-q_r}\rt) g$ follows from the projection formula
\begin{equation}
f_*(\beta)\wedge \alpha = f_*(\beta\wedge f^*\alpha)
\end{equation}
together with the formula
\begin{eqnarray}
f^*(\alpha_1\wedge\dots\wedge\alpha_r) = f^*(\alpha_1)\wedge \dots
\wedge f^*(\alpha_r) =  \\
\lt( \sum_{\sigma(1)=1}^r
n_{\sigma(1)}\pi_{\sigma(1)}^*(\alpha_1)\rt) \wedge \dots \wedge \lt(
\sum_{\sigma(2g-p)=1}^r
n_{\sigma(2g-p)}\pi_{\sigma(2g-p)}^*(\alpha_{2g-p})\rt). \nonumber
\end{eqnarray}

\end{proof}

\ps

\begin{thm}\label{sec-thmu}
The image $\Theta'$ is a divisor in $J(X)$ which is algebraically 
equivalent to $3\Theta$ and $\psi':\text{Sym}^2(F)\rightarrow \Theta'$
is birational.  The singular locus of $\Theta'$ has codimension
$2$.  For general $a\in J(X)$, we have that $\Theta'\cap(a+\Theta)$
is irreducible.
\end{thm}

\begin{proof}
For divisors on $J(X)$, algebraic equivalence is equivalent to
homological equivalence, so we shall establish that $\Theta'$ is
homologically equivalent to $3\Theta$.  

\

Using the notation of the previous lemma,
both $\Theta$ and $\Theta'$ are images of the subvariety $u(F)\times
u(F) \subset J \times J$ under the morphism $f_{(1,-1)}$ and
$f_{(1,1)}$ respectively.  By lemma~\ref{lem-cohom} we know there is a
linear map $g:H^4(J)\otimes H^4(J)\rightarrow H^2(J)$ such that 
\begin{eqnarray}
f_{(1,-1)*}(P.D.[u(F)\times u(F)]) = (1)^6(-1)^6 g(P.D.[u(F)\times
u(F)]), \\
f_{(1,1)*}(P.D.[u(F)\times u(F)]) = (1)^6(1)^6 g(P.D.[u(F)\times
u(f)]) \nonumber 
\end{eqnarray}
where $P.D.[u(F)\times u(F)]$ is the Poincar\'e dual of the homology
class of $u(F)\times u(F)$.  We know the mapping $\psi: F\times F
\rightarrow \Theta$ has degree $6$.  So if $\psi': F\times
F\rightarrow \Theta'$ has degree $d$, then it follows that
\begin{equation}
d P.D.[\Theta'] = f_{(1,1)*}(P.D.[u(F)\times u(F)]) =
f_{(1,-1)*}(P.D.[u(F)\times u(F)]) = 6 P.D.[\Theta].
\end{equation}
In other words, the cohomology class of $\Theta'$ is equal to
$\frac{6}{d}$ times the cohomology class of $\Theta$.  Of course $d$
is divisible by $2$ because $\psi'$ factors through $F\times F
\rightarrow \text{Sym}^2(F)$.  It remains to prove that $d$ is
precisely $2$.

\ps

Now let $\mathcal{G}':(\Theta')^{ns}\rightarrow \PP^{4\vee}$ be the
Gauss map defined on the nonsingular locus of $\Theta'$.  Of course
the differentials $d\psi$ and $d\psi'$ are simply given by $d u\circ
d\pi_1 - d u\circ d\pi_2$ and by $d u\circ d\pi_1 + d u\circ d\pi_2$
respectively.  In particular for each point $(L,M)\in F\times F$, the
subspace $\text{image}(d\psi)\subset T_0 J$ and
$\text{image}(d\psi)\subset T_0 J$ are equal, i.e. the composites
\begin{eqnarray}
(F\times F-\Delta)\xrightarrow{\psi} \Theta \xrightarrow{\mathcal{G}}
\PP^{4\vee} \\
F\times F\xrightarrow{\psi'} \Theta'
\xrightarrow{{\mathcal{G}}'} \PP^{4\vee}
\end{eqnarray}
are equal as rational maps; so both are equal to the rational map
$\Phi$.  Therefore, for each point $\lt([L_1],[L_2]\rt)$ whose image
under $\psi'$ lies in the nonsingular locus of $\Theta'$, we see that
the fiber of $\psi'$ containing $\lt([L_1],[L_2]\rt)$ is contained in
the fiber of $\Phi$ containing $\lt([L_1],[L_2]\rt)$.  So, we are
reduced to showing that, for generic $\lt([L_1],[L_2]\rt)$, if
$\psi'\lt([L_1],[L_2]\rt)=\psi'\lt([L_3],[L_4]\rt)$ with $L_1,L_2,L_3$
and $L_4$ contained in a smooth hyperplane section of $X$, then either
$\lt([L_1],[L_2]\rt)=\lt([L_3],[L_4]\rt)$ or
$\lt([L_2],[L_1]\rt)=\lt([L_4],[L_3]\rt)$.

\ps

A bit more generally, suppose that $\lt([L_1],[L_2]\rt)\in U$ is a pair
of skew lines.  We will show
that the fiber $\lt(\psi'\rt)^{-1}\lt(\psi'\lt([L_1],[L_2]\rt)\rt) = \lt\{
\lt([L_1],[L_2]\rt),\lt([L_2],[L_1]\rt)\rt\}$.  Indeed, suppose that 
$\psi'\lt([L_1],[L_2]\rt) =\psi'\lt([L_3],[L_4]\rt)$ and suppose
that $\lt([L_3],[L_4]\rt)\neq \lt([L_2],[L_1]\rt)$.  Then
$\psi\lt( [L_1],[L_4]\rt)=\psi\lt([L_3],[L_2]\rt)$.  Therefore
$\Phi\lt([L_1],[L_4]\rt)=\Phi\lt([L_3],[L_2]\rt)$.  Let us call this
common hyperplane $H$.  Then $L_1\subset H$ and $L_2\subset H$.  Therefore
$H=\Phi\lt([L_1],[L_2]\rt)$.  Since we have $L_4\neq L_1$ and $L_3\neq L_2$,
we conclude that $ \lt([L_1],[L_4]\rt), \lt([L_3],[L_2]\rt)\in U$.
Therefore, by theorem~\ref{thm-cg} (7), we have either $\lt([L_1],[L_4]\rt)=
 \lt([L_3],[L_2]\rt)$, or else there exists a line $l\in H$ such that
$l\cup L_1\cup L_2$ is the intersection of $X$ with a $\PP^2$.  But $L_1$
and $L_2$ cannot lie in a common $\PP^2$ since they are skew.  Therefore
we conclude that $\lt([L_1],[L_4]\rt)=
 \lt([L_3],[L_2]\rt)$, i.e. $\lt([L_3],[L_4]\rt)=\lt([L_1],[L_2]\rt)$.

\ps

So we deduce that $d=2$ and thus $\Theta'$ is algebraically equivalent 
to $3\Theta$.  But we deduce even more.  The image of 
$U-I\cap U$ in $\text{Sym}^2 F$ is smooth, let's call it $U'$.  Since
the map $\psi':U'\rightarrow\psi'(U')$ is bijective and since
$\text{rank}(d\psi') = \text{rank}(d\psi) = 4$ on $U - I\cap U$, 
it follows from Zariski's main
theorem~\cite[p. 288-289]{M} that 
$\psi'(U')$ is 
smooth and $U'\xrightarrow{\psi'}\psi'{U'}$ is an isomorphism.  Now the
complement of $U$ in $F\times F$ has codimension $2$.  Therefore
$\psi'(F\times F-U)$ has codimension at least $2$ inside of $\Theta'$.
But in fact $\psi'(I)$ also has codimension $2$.  Consider the rational
map $I\xrightarrow{\rho} F$ defined by sending a pair of incident lines
$\lt([L_1],[L_2]\rt)$ to the residual line $l$ such that
$\text{span}(L_1,L_2)\cap X = l\cup L_1\cup L_2$.  By the residuation trick
we have that the restriction of $\psi'$ to $I$ equals the composition 
of $\rho$ with the pointwise negative Abel-Jacobi map $-i$.  In
particular $\psi'(I)$ is just $-u(F)$ up to translation.  So $\psi'(I)$ has
codimension $2$ inside of $\Theta'$.  Therefore $\psi'(U')\subset\Theta'$
has complement of codimension $2$.  So the singular locus
of  $\Theta'$, has codimension $2$.

\ps

Since
$\Theta'\cap(a+\Theta)$ is a positive dimensional intersection of
ample divisors,  
$\Theta'\cap(a+\Theta)$ is connected.
Since $\Theta'\cap(a+\Theta)$ is a
complete intersection of Cartier divisors, it follows from
Hartshorne's connectedness theorem \cite[thereom 18.12]{E} that
$\Theta'\cap(a+\Theta)$ is connected in codimension $2$.  Thus to
prove that $\Theta'\cap (a+\Theta)$ is irreducible, it suffices to
prove that the singular locus of $\Theta'\cap (a+\Theta)$ has
codimension at least $2$.
By the Bertini-Kleiman theorem~\cite[theorem III.10.8]{H}, we
have that for general $a\in J(X)$, the intersection $\Theta'\cap(a+\Theta)$
is smooth away from the intersection of each divisor with the singular 
locus of the other divisor.  But by the last paragraph and by
theorem~\ref{thm-tj}, we see that the singular loci of $\Theta$ and
$\Theta'$ both have codimension at least $2$ (in $\Theta$ or $\Theta'$ 
respectively).  So it follows that for $a$ general, the singular locus 
of $\Theta'\cap(a+\Theta)$ has codimension at least $2$.  Thus we
conclude that for general $a\in J(X)$, the intersection
$\Theta'\cap(a+\Theta)$ is irreducible.

\end{proof}

We also use the following enumerative lemma, which is proved
in~\cite[lemma 4.2]{HRS1}.  

\ps

\begin{lem}\label{lem-enum}  Suppose that $C\subset X$ is a smooth
curve of genus $g$ 
and degree $d$.  Let $B_C\subset F$ denote the scheme parametrizing
lines in $X$ which intersect $C$ in a scheme of degree $2$ or more.
Define $b(C)=\frac{5d(d-3)}{2}+6-6g$.  If
$B_C$ is not positive dimensional and if $b(C)\geq 0$, then the degree of
$B_C$ is $b(C)$.  
\end{lem}

\section{Twisted cubics}~\label{sec-30}

We  now begin in earnest the analysis of the geometry of cubic,
quartic and quintic curves on our cubic threefold 
$X \subset \P^4$. In each case our goal is the proof of the Main
Theorem \ref{mainthm} for the variety 
$\Hdg{d}{g}(X)$ parametrizing curves of this degree and genus. We
start with the variety $\Hdg{3}{0}(X)$ 
parametrizing rational curves of degree 3---that is \emph{twisted cubics}. 

\ps

In ~\cite[theorem 4.4]{HRS1} we prove the following result

\begin{thm}\label{thm-irr30} 
The space $\Hdg{3}{0}(X)$ is a smooth,
irreducible $6$-dimensional variety.
\end{thm}

We have a morphism 
\begin{equation}
\Hdg{3}{0}(\PP^4)\xrightarrow{\sigma^{3,0}} \PP^{4\vee}
\end{equation}
defined by sending $[C]$ to $\text{span}(C)$.
This morphism makes $\Hdg{3}{0}(\PP^4)$ 
into a locally trivial bundle over $\PP^{4\vee}$ with fiber
$\Hdg{3}{0}(\PP^3)$ .  Recall from section~\ref{sec-123} that we defined
$X^\vee
\subset \PP^{4\vee}$ to be the dual variety of $X$ which parametrizes
tangent hyperplanes to $X$ and we defined $U$ to
be the complement of $X^\vee$ in $\PP^{4\vee}$.  Then we define 
$\Hdg{3}{0}_U(X)$ to be the
open subscheme of $\Hdg{3}{0}(X)$ which parametrizes twisted cubics,
$C$, in $X$ 
such that $\sigma^{3,0}([C])\in U$.  By the graph construction we may
consider $\Hdg{3}{0}_U(X)$ as a locally closed subvariety of
$\HI{3t+1}{X}\times U$.  
Let $\overline{\mathcal{H}}\subset \HI{3t+1}{X}\times U$ denote the
closure of $\Hdg{3}{0}_U(X)$ with the reduced induced scheme structure.
Denote by $\overline{\mathcal{H}} \xrightarrow{f} U(X)$ the projection map. 

\begin{thm}\label{thm-st30}  
Let
$\overline{\mathcal{H}}\xrightarrow{f''} U'\xrightarrow{f'} 
U(X)$ 
be the Stein factorization of $\overline{\mathcal{H}}\xrightarrow{f} U(X)$.
Then 
$\overline{\mathcal{H}}\xrightarrow{f''} U'$ is isomorphic to a $\PP^2$-bundle
$\PP_{U'}(E)\xrightarrow{\pi} U'$ with $E$ a locally free sheaf of rank $3$.
And $U'\xrightarrow{f'} U$ is an unramified finite morphism of degree
$72$.  Moreover, the Abel-Jacobi map $\overline{\mathcal{H}}\xrightarrow{i}
J(X)$ 
factors as $\overline{\mathcal{H}}\xrightarrow{f''} U'\xrightarrow{i'}
J(X)$ where 
$U'\xrightarrow{i'} J(X)$ is a birational morphism of $U'$ to a translate
of $\Theta$.
\end{thm} 

\begin{proof}
This theorem is  ~\cite[theorem 4.5]{HRS1}, and is proved there.  For
brevity's sake, we recall just the sketch of the proof.

\ps

One can form the subvariety $\mathcal{X}\subset U\times X$ as the
universal smooth hyperplane section of $X$.  Each fiber of
$\mathcal{X}\rightarrow U$ is a smooth cubic surface.  We associate to
this family of cubic surfaces the finite \'etale morphism
\begin{equation}
\rho: \text{Pic}^{3,0}(\mathcal{X}/U) \rightarrow U
\end{equation}
whose fiber over a point $[H]\in U$ is simply the set of divisor
classes $D$ on the smooth cubic surface $H\cap U$ such that $D.D =1$
and $D.h=3$
with $h$ the hyperplane class.  Now the Weyl group $\mathcal{W}(E_6)$
acts transitively on the set of such divisor classes $D$. It follows
that with respect to the model of a cubic surface as 
$\PP^2$ blown up at 6 points, each divisor class $D$ above corresponds
to the class of a line on $\PP^2$.  From this it follows that the
induced morphism 
\begin{equation}
g:\Hdg{3}{0}(X)\rightarrow \text{Pic}^{3,0}(\mathcal{X}/U))
\end{equation}
is surjective and isomorphic to an open subset of a $\PP^2$-bundle on
$\text{Pic}^{3,0}(\mathcal{X}/U)$.  In particular, since
$u:\Hdg{3}{0}(X)\rightarrow J(X)$ contracts all rational curves, there
is an induced morphism $u':\text{Pic}^{3,0}(\mathcal{X}/U)\rightarrow
J(X)$ such that $u=u'\circ g$.  

\ps

It only remains to show that $u'$ maps
$\text{Pic}^{3,0}(\mathcal{X}/U)$ birationally to its image.  This is
proved by analyzing ``Z''s of lines, i.e. configurations of lines on a
cubic surface $H\cap X$ whose dual graph is just the connected graph with 3
vertices and no loops.  Using ~\ref{thm-cg}, part 5, we conclude
that two ``Z''s have the same image in $J(X)$ iff they are in the same
linear equivalence class on $H\cap X$.  
\end{proof}

\begin{cor}\label{cor-aj30}  The Abel-Jacobi map
$u_{3,0}:\Hdg{3}{0}(X)\rightarrow J(X)$ dominates a translate of
$\Theta$ and is birational to a $\PP^2$-bundle over its image.
\end{cor}

\begin{proof}  Since $\Hdg{3}{0}(X)$ is irreducible, $\Hdg{3}{0}_U(X)$ 
is dense in $\Hdg{3}{0}(X)$.
\end{proof}

\subsection{Quartic Elliptic Curves}

Recall that the normalization of $\HI{2t+1}{X}$ is isomorphic to the
$\PP^2$-bundle $\PP 
Q\rightarrow F$ which parametrizes pairs $(L,P)$ which $L\subset X$
a line and $P\subset\PP^4$ a $2$-plane containing $L$.    
Let $A\xrightarrow{g} \PP Q$ denote the $\PP^1$-bundle
which parametrizes triples $(L,P,H)$ with $H$ a hyperplane containing
$P$.  Let $I_{4,1}\xrightarrow{h} A$ denote the $\PP^4$-bundle parametrizing
4-tuples $(L,P,H,Q)$ where $Q\subset H$ is a quadric surface
containing the conic $C\subset X\cap P$. 
Notice that $I_{4,1}$ is smooth
and connected of dimension $4+1+4=9$.

\ps

Let $D\subset I_{4,1}\times X$
denote the intersection of the universal quadric surface over $I_{4,1}$ with 
$I_{4,1}\times X\subset I_{4,1}\times \PP^4$.  Then $D$ is a local complete
intersection scheme.  By the Lefschetz hyperplane theorem, $X$
contains no quadric surfaces; therefore $D\rightarrow I_{4,1}$ has constant
fiber dimension $1$ and so is flat.  Let $D_1\subset I_{4,1}\times X$ denote
the pullback from $\PP Q\times X =\HI{2t+1}{X}\times X$ of the
universal family of conics.  Since $I_{4,1}\times X\rightarrow \PP Q\times
X$ is smooth and the universal family of conics is a local complete
intersection which is flat over $\PP Q$, we conclude that also $D_1$ is a
local complete intersection which is flat over $I_{4,1}$.  Clearly
$D_1\subset D$.  Thus by corollary~\cite[corollary 2.7]{HRS1}, we see that the
residual $D_2$ of $D_1\subset D$ is Cohen-Macaulay and flat over $I_{4,1}$.

\ps

By the base-change property in corollary~\cite[corollary 2.7]{HRS1},
we see that 
the fiber of $D_1\rightarrow I_{4,1}$ over a point $(L,P,H,Q)$ is simply the
residual of $C\subset Q\cap X$.  If we choose $Q$ to be a smooth
quadric, i.e. $Q\cong \PP^1\times \PP^1$, then $C\subset Q$ is a
divisor of type $(1,1)$ and $X\cap Q\subset Q$ is a divisor of type
$(3,3)$.  Thus the residual curve $E$ is a divisor of type $(2,2)$,
i.e. a quartic curve of arithmetic genus 1.  Thus $D_2\subset I_{4,1}\times
X$ is a family of connected, closed subschemes of $X$ with Hilbert
polynomial $4t$.  So we have an induced map $f:I_{4,1}\rightarrow
\HI{4t}{X}$.  

\ps

\begin{prop}\label{lem-qell}  The image of the morphism above
$f:I_{4,1}\rightarrow\HI{4t}{X}$ is the closure $\overline{\Hdg{4}{1}}(X)$ 
of $\Hdg{4}{1}(X)$.  Moreover the open set $f^{-1}\Hdg{4}{1}(X)\subset 
I_{4,1}$ is a $\PP^1$-bundle over $\Hdg{4}{1}(X)$.  Thus $\Hdg{4}{1}(X)$ is
smooth and connected of dimension 8.
\end{prop}

\begin{proof}
This is ~\cite[proposition 5.1]{HRS1}.
\end{proof}

If $(L,P,H,Q)$ is a point in the fiber over $[E]$, then
$H=\text{span}(E)$.  And since there is no rational curve in $F$, we
also have that $L$ is constant in the fiber.
So we have a well-defined morphism $m:\Hdg{4}{1}(X)\rightarrow \PP
Q^\vee$, where $\PP Q^\vee$ is the $\PP^2$-bundle over $F$ parametrizing pairs
$([L],[H]), L\subset H$.  For a general $H$, the intersection $Y=H\cap X$ is
a smooth cubic surface.  And the fiber $m^{-1}([L],[H])$ is an open
subset of the complete linear series $|\OO_Y{L+h}|$, where $h$ is the
hyperplane class on $Y$.  Thus $m:\Hdg{4}{1}(X)\rightarrow \PP Q^\vee$
is a morphism of smooth connected varieties which is birational to a
$\PP^4$-bundle.  Composing $m$ with the projection $\PP Q^\vee$ yields
a morphism $n:\Hdg{4}{1}(X)\rightarrow F$ which is birational to a
$\PP^4$-bundle over a $\PP^2$-bundle.  

\begin{cor}~\label{cor-qell}
By the residuation trick, we conclude 
the Abel-Jacobi map
$u_{4,1}:\Hdg{4}{1}(X)\rightarrow J(X)$ is equal, up to a
fixed translation, to the composite:
\begin{equation}\begin{CD}
\Hdg{4}{1}(X) @> n >> F @> u_{1,0} >> J(X).
\end{CD}\end{equation}
Thus the general fiber of $u_{4,1}$ equals the general fiber
of $n$, and so is isomorphic to an open subset of a $\PP^4$-bundle
over $\PP^2$.  
\end{cor}

\section{Cubic scrolls and applications}~\label{sec-cuscr}

\subsection{Preliminaries on cubic scrolls}~\label{subsec-prelim}

In the next few sections we will use residuation in a cubic surface
scroll. We start by collecting 
some basic facts about these surfaces.

\ps

There are several equivalent descriptions of cubic scrolls.  
\begin{enumerate}
\item A cubic scroll $\Sigma \subset \PP^4$ is a connected, smooth
surface with Hilbert polynomial $P(t) = \frac{3}{2}t^2 + \frac{5}{2}t
+ 1$. 
\item A cubic scroll $\Sigma\subset
\PP^4$ is the determinantal variety defined by the $2\times 2$ minors of 
a matrix of linear forms:
\begin{equation}\lt[ \begin{array} {ccc}
 L_1 & L_2 & L_3 \\
 M_1 & M_2 & M_3 
\end{array} \rt]\end{equation}
such that for each row or column, the linear forms in
that row or column are linearly independent 
\item A cubic scroll $\Sigma\subset \PP^4$ is the \emph{join} of an
isomorphism $\phi:L\rightarrow C$.  Here $L\subset \PP^4$ is a line and
$C\subset \PP^4$ such that $L\cap
\text{span}(C) =\emptyset$.  The join of $\phi$ is defined as the
union over all $p\in L$ of the line $\text{span}(p,\phi(p))$.  
\item A cubic scroll $\Sigma\subset \PP^4$ is the image of a morphism $f:\PP
E\rightarrow \PP^4$ where $E$ is the rank 2 vector bundle on $\PP^1$, 
$E=\OO_{\PP^1}(-1)\oplus \OO_{\PP^1}(-2)$, and the morphism $f:\PP
E\rightarrow \PP^4$ is such that $f^*\OO_{\PP^4}(1) = \OO_{\PP E}(1)$
and the pullback map $H^0(\PP^4,\OO_{\PP^4}(1))\rightarrow H^0(\PP
E,\OO_{\PP E}(1)$ is an isomorphism.
\item A cubic scroll $\Sigma\subset \PP^4$ is as
a minimal variety, i.e. $\Sigma\subset \PP^4$ is any smooth connected 
surface with $\text{span}(\Sigma)=\PP^4$ which has the minimal
possible degree for such a surface, namely $\text{deg}(\Sigma)=3$. 
\item A cubic scroll $\Sigma\subset \PP^4$ is a smooth surface residual to a
2-plane $\Pi$ in the base locus of a pencil of quadric hypersurfaces
which contain $\Pi$.  
\end{enumerate}

\ps

From the fourth description $\Sigma = f(\PP E)$ we see that
$\text{Pic}(\Sigma) = \text{Pic}(\PP E) \cong \ZZ^2$.  Let $\pi:\PP
E\rightarrow \PP^1$ denote the projection morphism and let
$\sigma:\PP^1\rightarrow \PP E$ denote the unique section whose image
$D=\sigma(\PP^1)$ has self-intersection $D.D=-1$.  Then $f(D)$ is a
line on $\Sigma$ called the \emph{directrix}.  And for each $t\in
\PP^1$, $f(\pi^{-1}(t))$ is a line called a \emph{line of the ruling} of
$\Sigma$.  Denote by $F$ the divisor class of any $\pi^{-1}(t)$.  Then
$\text{Pic}(\Sigma) = \ZZ\{D,F\}$ and the intersection pairing on
$\Sigma$ is determined by $D.D = -1, D.F =1, F.F =0$.  The hyperplane
class is $H=D+2F$ and the canonical class is $K = -2D-3F$.  

\ps

Using the fourth description of a cubic scroll, we see that any two
cubic scrolls differ only by the choice of the isomorphism
$H^0(\PP^4,\OO_{\PP^4}(1)) \rightarrow H^0(\PP E,\OO_{\PP E}(1))$.
Therefore any two cubic scrolls are
conjugate under the action of $\text{PGL}(5)$.  So the open set
$U\subset \text{Hilb}_{P(t)}(\PP^4)$ parametrizing cubic scrolls is a
homogeneous space for 
$\text{PGL}(5)$, in particular it is smooth, connected and rational.
So the Abel-Jacobi map $U\rightarrow J(X)$ associated to the family of
intersections $\Sigma\cap X\subset X$ is a constant map.  

\subsection{Cubic Scrolls and Quartic Rational
Curves}~\label{subsec-scr40}

Recall that $\text{Pic}(\Sigma) = \ZZ\{D,F\}$ where $D$ is the
directrix and $F$ is the class of a line of ruling.  The intersection
product is given by $D^2 = -1, D.F = 1, F^2 =0$.  The canonical
class is given by $K_\Sigma = -2D-3F$ and the hyperplane class is
given by $H=D+2F$.  The linear system $|F|$ is nef because it is the
pullback of $\OO_{\PP^1}(1)$ under the projection
$\pi:\Sigma\rightarrow \PP^1$.  Similarly, $|D+F|$ is nef because it
contains all the conics obtained as the residuals to lines of the
ruling in $|H|$.  Thus for any effective curve class $aD+bF$ we have
the two inequalities $a=(aD+bF).F \geq 0$, $b=(aD+bF).(D+F) \geq 0$.  

\ps

Suppose that $C\subset \Sigma$ is an effective
divisor of degree $4$ and arithmetic genus $0$.  By the adjunction
formula 
\begin{equation}
K_\Sigma.[C] + [C].[C] = 2p_a - 2 = -2.
\end{equation}
So if $[C]=aD+bF$, then we have the conditions
\begin{equation}
a\geq 0, b\geq 0, a+b = 4, a^2 -2ab +a+2b = 2.
\end{equation}
It is easy to check that there are precisely two solutions $[C]=2D+2F,
[C]=D+3F$.  We will see that both possibilities occur and describe
some constructions related to each possibility.

\ps

\begin{lem}\label{lem-g1240}  Let $C\subset \PP^4$ be a smooth quartic 
rational curve and let $V\subset |\OO_C(2)|$ be a pencil of degree
$2$-divisors on $C$ without basepoints.  There exists a unique map of a
Hirzebruch surface $\mathbf{F}_1$, $f:\Sigma\rightarrow \PP^4$ such
that $(f^* H)^2=3$  
and a factorization $i:C\rightarrow \Sigma$ of $C\rightarrow \PP^4$
such that $i(C)\sim 2D+2F$ and the pencil of degree $2$ divisors
$F\cap C$ is the pencil $V$.
\end{lem}

\begin{proof}
This is ~\cite[lemma 6.7]{HRS1}.  
\end{proof}

\textbf{Remark}
While we are at it, let's mention a specialization of the construction
above, namely what happens if $V$ is not basepoint free.  Then
$V=p+|\OO_C(1)|$ where $p\in C$ is some basepoint.  Consider
the projection morphism $f:\PP^4\rightarrow \PP^3$ obtained by
projection from $p$ (this is a rational map undefined at $p$).  The
image of $C$ is a rational cubic curve $B$ (possibly a singular plane
cubic).  Consider the cone $\Sigma'$ in $\PP^4$ over $B$ with vertex
$p$.  This 
surface contains $C$.  If we blowup $\PP^4$ at $p$, then the proper
transform of $\Sigma'$ in $\widetilde{\PP^4}$ is a surface whose
normalization $\Sigma$ is a 
Hirzebruch surface $\mathbf{F}_3$ (normalization is only necessary if
$B$ is a plane curve).  The directrix $D$ of $\Sigma$ is the pullback of
the exceptional divisor of $\widetilde{\PP^4}$.  The inclusion
$C\subset \Sigma'$ induces a factorization $i:C\rightarrow \Sigma$ of
$C\rightarrow \PP^4$, and $[i(C)] = D + 4F$.  The intersection of $D$
and $i(C)$ is precisely the point $p$.  And the linear system $i^*|F|$
is exactly $|\OO_C(1)|$.  

\ps

Next we consider the case of a rational curve $C\subset \Sigma$ such
that $[C]=D+3F$.  

\begin{lem}~\label{lem-scr40}  
Let $C\subset \PP^4$ be a smooth quartic rational curve 
and let $L\subset \PP^4$ be a line such that $L\cap C=Z$ is a degree 2
divisor.  
Let $\phi:C\rightarrow L$ be an isomorphism such that $\phi(Z)=Z$ and
$\phi|_Z$ is the identity map.  Then there exists a unique triple
$(h,i,j)$ where $h:\Sigma\rightarrow \PP^4$ is a finite map of a 
Hirzebruch surface $\Sigma~\equiv~\mathbf{F}_1$, and 
$i:C\rightarrow\PP^4$, $j:L\rightarrow\PP^4$ are factorizations of
$C\rightarrow \PP^4,L\rightarrow \PP^4$
such that $j(L)=D$ is the directrix, such that $[i(C)]=D+3F$ and such
that the composition of $i:C\rightarrow \Sigma$ with the projection
$\pi:\Sigma\rightarrow D$ equals $j\circ \phi$.  
\end{lem}

\begin{proof}
This is ~\cite[lemma 6.8]{HRS1}.
\end{proof}

\subsection{Cubic Scrolls and Quintic Elliptics}~\label{subsec-scr51}

Recall that our fourth description of a cubic scroll was the image
of a morphism $f:\Sigma \rightarrow \PP^4$ where $\Sigma$ is the
Hirzebruch surface $\mathbb{F}_1$, $f^*\OO(1)\sim \OO_{\Sigma}(1) =
\OO_{\PP E}(D+2F)$, and $f:\Sigma \rightarrow \PP^4$ is given by the
complete linear series of $\OO_{\Sigma}(D+2F)$.  In the next 
sections it
will be useful to weaken this last condition.

\ps

\begin{defn} A \emph{cubic scroll} in $\PP^n$ is a finite morphism
$f:\Sigma\rightarrow \PP^n$ where $\Sigma$ is isomorphic to the
Hirzebruch surface $\mathbb{F}_1$ and such that $f^*\OO_{\PP^n}(1)$ is
isomorphic to $\OO_{\Sigma}(D+2F)$. 
\end{defn}

Let $H = D+2F$ denote the pullback of the hyperplane class.
Now suppose that $E\subset \Sigma$ is an effective Cartier divisor
with $p_a(E)=1$ and $E.H = 5$.  Since $F$ and $D+F$ are effective and
move, we have $E.F, E.(D+F) \geq 0$.  Writing $E=aD + bF$ we see
$(a,b)$ satisfies the relations $a,b \geq 0, a+b = 5$ and
$a(b-3)+b(a-2) -a(a-2) = 0$.  These relations give the unique solution
$E=2D+3F = -K$.  In particular, if $E$ is smooth then
$\pi:E\rightarrow \PP^1$ is a finite morphism of degree 2, i.e. a
$g^1_2$ on $E$.  Thus a pair $(f:\Sigma\rightarrow \PP^n, E\subset
\Sigma)$ of a cubic scroll and a quintic elliptic determines a pair
$(g:E\rightarrow \PP^n, \pi:E\rightarrow \PP^1)$ where $g:E\rightarrow
\PP^n$ is a quintic elliptic and $\pi:E\rightarrow \PP^1$ is a degree
2 morphism.  

\ps

Suppose we start with a pair $(g:E\rightarrow \PP^n, \pi:E\rightarrow
\PP^1)$ where $g:E\rightarrow \PP^n$ is an embedding of a quintic
elliptic curve and $\pi:E\rightarrow \PP^1$ is a degree 2 morphism.
Consider the rank 2 vector bundle $\pi_* g^* \OO_{\PP^n}(1)$.  

\ps

\begin{lem} Suppose $E$ is an elliptic curve and $\pi:E\rightarrow
\PP^1$ is a degree 2 morphism.  Suppose $L$ is an invertible sheaf on
$E$ of degree $d$.  Then we have

\begin{equation}
\pi_* L \cong \lt\{ \begin{array}{ll}
 \OO_{\PP^1}(e) \oplus \OO_{\PP^1}(e-1) & d = 2e+1, \\
 \OO_{\PP^1}(e) \oplus \OO_{\PP^1}(e-2) & d=2e, L\cong
        \pi^*\OO_{\PP^1}(e), \\
 \OO_{\PP^1}(e-1) \oplus \OO_{\PP^1}(e-1) & d=2e, L\not\cong \pi^*
        \OO_{\PP^1}(e) 
\end{array} \rt.
\end{equation}
\end{lem}

\begin{proof}
This is ~\cite[lemma 6.10]{HRS1}.
\end{proof}

By the lemma we see that the vector bundle $G:=\pi_* g^*
\OO_{\PP^n}(1)$ is isomorphic to $\OO_{\PP^1}(1)~\oplus~
\OO_{\PP^1}(2)$.  Associated to the linear series
$\OO_E^{n+1}\rightarrow g^*\OO_{\PP^n}(1)$ defining the embedding $g$,
we have the push-forward linear series $\OO_{\PP^1}^{n+1}\rightarrow
G$.  Since $g$ is an embedding, for each pair of points $\{p,q\}
\subset E$ (possibly infinitely near), we have that
$\OO_E^{n+1}\rightarrow g^*\OO_{\PP^n}(1)|_{\{p,q\}}$ is surjective.
In particular taking $\{p,q\}=\pi^{-1}(t)$ for $t\in \PP^1$, we
conclude that $\OO_{\PP^1}^{n+1}\rightarrow F|_t$ is surjective.  Thus
we have an induced morphism $\PP G^\vee \rightarrow \PP^n$ which pulls
back $\OO_{\PP^n}(1)$ to $\OO_{\PP G^\vee }(1)$.  Let us
denote $\Sigma := \PP G^\vee $ and let us denote the morphism by
$f:\Sigma\rightarrow \PP^n$.  Abstractly $\Sigma$ is isomorphic to
$\mathbb{F}_1$ and $f:\Sigma\rightarrow \PP^n$ is a cubic scroll.

\ps

The tautological map $\pi^*\pi_* g^*\OO_{\PP^n}(1)\rightarrow
g^*\OO_{\PP^n}(1)$ is clearly surjective.  Thus there is an induced
morphism $h:E\rightarrow \Sigma$.  Chasing definitions, we see that
$g=f\circ h$.  So we conclude that given a pair $(g:E\rightarrow
\PP^n, \pi:E\rightarrow \PP^1)$ as above, we obtain a pair
$(f:\Sigma\rightarrow \PP^n, h:E\rightarrow \Sigma)$.  Thus we
have prove the following:

\begin{lem}\label{lem-equiv1}
There is an equivalence between the collection of pairs
$(f:\Sigma\rightarrow \PP^n,h:E\rightarrow \Sigma)$ with
$f:\Sigma\rightarrow \PP^n$ a cubic scroll and $f\circ h:E\rightarrow
\PP^n$ an embedded quintic elliptic curve and the collection of pairs
$(g:E\rightarrow \PP^n, \pi:E\rightarrow \PP^1)$ where $g:E\rightarrow
\PP^n$ is an embedded quintic elliptic curve and $\pi:E\rightarrow
\PP^1$ is a degree 2 morphism.  
\end{lem}

Stated more precisely, this gives an
isomorphism of the parameter schemes of such pairs, but we don't need
such a precise result.

\subsection{Cubic Scrolls and Quintic Rational
Curves}~\label{subsec-scr50}

If one carries out the analogous computations as at the beginning of
subsection~\ref{subsec-scr40} one sees that the only effective divisor
classes $aD+bF$ on a cubic scroll $\Sigma$ with degree $5$ and
arithmetic genus $0$ are $D+4F$ and $3D+2F$.  But the divisor class
$3D+2F$ cannot be the divisor of an irreducible curve because
$(3D+2F).D = -1$.  Thus if $C\subset \Sigma$ is an irreducible curve
of degree $5$ and arithmetic genus $0$, then $[C]=D+4F$.  

\begin{lem}~\label{lem-scr50}
Let $C\subset \PP^4$ be a smooth quintic rational curve
and let $L\subset \PP^4$ be a line such that $L\cap C$ is a degree 3
divisor $Z$.  Let $\phi:C\rightarrow L$ be an isomorphism such that
$\phi(Z)=Z$ and $\phi|_Z$ is the identity map.  Then there exists a
unique triple $(h,i,j)$ such that $h:\Sigma\rightarrow\PP^4$ is a
finite map of a Hirzebruch surface $\Sigma~\equiv~\PP^4$, and
$i:C\rightarrow\PP^4, j:L\rightarrow\PP^4$ are factorizations of
$C\rightarrow\PP^4, L\rightarrow\PP^4$ such that $j(L)=D$ is the
directrix, such that $[i(C)]=D+4F$ and such that the composition of
$i$ with the projection $\pi:\Sigma \rightarrow D$ equals $j\circ
\phi$.  
\end{lem}

\begin{proof}
This is ~\cite[lemma 6.12]{HRS1}.
\end{proof}

\section{Quartic Rational Curves}

In this section we will prove that the Abel-Jacobi map
$u_{4,0}:\Hdg{4}{0}(X)\rightarrow J(X)$ is dominant and the general
fiber is irreducible of dimension $3$.  In a later section we will
prove that the general fiber is unirational.

\ps

The construction we use to understand quartic rational curves is as
follows.  For any quartic rational curve
$C\subset X$, define $A_C\subset F$ to be the scheme parametrizing
$2$-secant lines to $C$.  By lemma~\ref{lem-enum}, $A_C$ is either
positive-dimensional or else has length $16$.  Define $I\subset
\Hdg{4}{0}(X)\times \HI{2}{F}$ to be the space of pairs $([C],[Z])$
where $Z\subset A_C$ is a $0$-dimensional length 2 subscheme of $A_C$.  
Denote by $\Lambda\subset Z\times \PP^4$ the flat family of lines
determined by $Z\subset F$.  

\begin{defn} We say $Z$ is \emph{planar} if there exists a $2$-plane
$P\subset \PP^4$ such that $\Lambda\subset Z\times P$.  We say $Z$ is
\emph{nonplanar} if $Z$ is not planar.
\end{defn}

By lemma 7.1, lemma 7.2 and theorem 7.3 of ~\cite{HRS1}, we have
the following result

\begin{thm}~\label{thm-irr40} 
The morphism $I\rightarrow \Hdg{4}{0}(X)$ is generically
finite, $I$ is irreducible of dimension $8$ and therefore
$\Hdg{4}{0}(X)$ is irreducible of dimension $8$.  For a general pair
$([C],[Z])\in I$, $C$ is nondegenerate and $Z$ is reduced and nonplanar.
\end{thm}

Let $([C],[Z])\in I$ be a general pair.
Now we may also consider $Z$ as a subscheme of
$\text{Sym}^2(C)\cong\PP^2$.  Since $Z$ is non-planar, in particular
the span of $Z\subset \PP^2$ yields a pencil of degree $2$ divisors on
$C$ without basepoints.  By lemma~\ref{lem-g1240}, there is a cubic
scroll $f:\Sigma\rightarrow\PP^4$ along with a factorization
$i:C\rightarrow\PP^4$ such that the pencil of degree $2$ divisors on
$C$ is just the pencil of intersections of $C$ with the lines of
ruling of $\Sigma$.  The residual of $C\cup \Lambda$ in $\Sigma\cap X$
is a curve in $\Sigma$ of degree $3$ and arithmetic genus $0$.  As a
corollary of the proof of ~\cite[theorem 7.3]{HRS1}, for $([C],[Z])$
a general pair, the residual curve $D$ is a twisted cubic curve.
As in theorem~\ref{thm-st30}, let $U'\subset J(X)$ denote the
Abel-Jacobi image of the locus of twisted cubics $D\subset X$ such
that $\text{span}(D)\cap X$ is a smooth cubic surface.  And recall
from theorem~\ref{sec-thmu} that the Abel-Jacobi map
$\HI{2}{F}\rightarrow \Theta'$ is birational.  Define $R=\Theta'\times
U'$ and define a rational transformation 
\begin{equation}
h:I\rightarrow R, \lt([C],[Z]\rt)\mapsto (u[Z],u([D])).
\end{equation}

\begin{lem}  The rational transformation $h:I\rightarrow R$ is
birational.
\end{lem}

\begin{proof}  By Zariski's main theorem, it suffices to prove that
for the general element $(x,y)\in R$, there is a unique pair
$([C],[Z])$ such that $h([C],[Z])=(x,y)$.  Now $x=\{L_1,L_2\}$ is a
general pair of disjoint lines in $X$.  And $y$ is a linear
equivalence class $|D|$ of twisted cubics on a general hyperplane
section $H\cap X$.  Then $L_i\cap H=\{p_i\}$ is a general point on $H\cap X$
(since every point of $X$ lies on a line, we may assume this point is
general on $H\cap X$).  And there is a unique twisted cubic $D\subset
H\cap X$ in the linear equivalence class $|D|$ and which contains the
two points $p_1$ and $p_2$.  If $h([C],[Z])=(x,y)$, then the residual
to $C\cup L_1\cup L_2$ in $\Sigma$ can only be $D$.  

\ps

Now consider the hyperplane $H'=\text{span}(L_1,L_2)$.  The
intersection $H'\cap D$ consists of three points $p_1,p_2$ and a third
point $q$.  Moreover the directrix $M$ of $\Sigma$ is contained in
$H$.  And by a divisor class calculation, $M\cap D$ consists of a
point other than $p_1, p_2$.  The only possibility is that $M\cap D =
\{q\}$.  Notice that given two skew lines $L_1, L_2$ in a $3$-plane
$H'$ and given a point $q\in H'$ not lying on $L_1\cup L_2$, there is
a unique line $M\subset H'$ which contains $q$ and intersects each of
$L_1,L_2$.  Indeed, if we project from $q$ then $L_1,L_2$ project to
distinct lines in $H/q\cong \PP^2$ which intersect in a unique point.
And $M$ is simply the cone over this point.  So we conclude that the
directrix line $M$ is uniquely determined by $(x,y)$.  

\

Finally, projection $\Sigma\rightarrow M$ to the directrix determines
an isomorphism 
$\phi:D\rightarrow M$ such that for each $r\in D,
\text{span}(r,\phi(r))$ is a line of the ruling of $\Sigma$, in
particular $\phi(q)=q$.  Conversely, given a line $M$ which intersects
$C$ in one point $q$ and given an isomorphism $\phi:C\rightarrow M$
such that $\phi(q)=q$, then the union of the lines
$\text{span}(r,\phi(r))$ is a cubic scroll $\Sigma$.  Thus to uniquely
specify the cubic scroll $\Sigma$, we have only to determine $\phi$.
But notice also that we already know $\phi(p_1),\phi(p_2)$ are the
unique points of intersection of $M$ with $L_1$ and $L_2$
respectively.  Thus $\phi$ is also uniquely determined by $(x,y)$.
Altogether we conclude that $\Sigma$ is uniquely determined by
$(x,y)$.  But then we can recover/construct $C$ as the residual to
$D\cup L_1\cup L_2$ in $\Sigma\cap X$.  This proves that $h$ is
birational.
\end{proof}

\ps

\begin{thm}\label{thm-irrfib} 
The composite of $I\rightarrow
\Hdg{4}{0}(X)$ with the Abel-Jacobi map
$u_{4,0}:\Hdg{4}{0}(X)\rightarrow J(X)$ is dominant and the
general fiber is irreducible.  Therefore $u_{4,0}$ is dominant and the 
general fiber is irreducible.
\end{thm}

\begin{proof}
By the residuation trick, the composite equals (as a rational map) the
pointwise inverse of the composition
\begin{equation}
I\xrightarrow{h} R = \Theta'\times U' \rightarrow J(X)
\end{equation}
where the second map is the restriction to $\Theta'\times U'\subset
J(X)\times J(X)$ of the addition map.  Clearly $\Theta'+U'= J(X)$ since
$\Theta'$ is a divisor. $U'\subset \Theta$ is a Zariski-dense open
set, and $\Theta$ is not contained in any translate of
$\Theta'$.  Thus $I\rightarrow J(X)$ is dominant.  

\ps

Using the fact that $\Theta$ is a symmetric divisor, we see that the
general fiber of the map $\Theta'\times\Theta \rightarrow J(X)$ is an
intersection $\Theta\cap (a+\Theta')$.  By theorem~\ref{sec-thmu}, for 
general $a$ this intersection is irreducible.  Thus we conclude that
the general fiber of $I\rightarrow J(X)$ is irreducible.
This proves the theorem.
\end{proof}

\section{Quintic Elliptics}~\label{subsec-irr51}

In this section we will prove that the Abel-Jacobi map
$u_{5,1}:\Hdg{5}{1}(X)\rightarrow J(X)$ is dominant and the general
fiber is an irreducible $5$-fold.  In the next section we will see
that the fibers are unirational.

\ps

The construction we use to understand quintic elliptics is as
follows.  Define $g:\widetilde{H}\rightarrow \Hdg{5}{1}(X)$ to be the
relative $\text{Pic}^2$ of the universal family of elliptic curves
$\mathcal{C}\rightarrow \Hdg{5}{1}(X)$.  By lemma~\ref{lem-equiv1}
$\widetilde{H}$ is also the parameter space for pairs $(f,h)$ where
$f:\Sigma\rightarrow \PP^4$ is a generalized cubic scroll and
$h:C\rightarrow \Sigma$ is a curve such that $f(h(C))\subset X$ is a
smooth quintic elliptic.  The residual to $C$ in $\Sigma\cap X$ is a
curve $C'$ of degree $4$ and arithemtic genus $0$.  

\begin{thm}~\label{thm-irr51}  
The scheme $\widetilde{H}$ is irreducible of dimension
$11$.  For a general pair $(f,h)\in \widetilde{H}$, the residual curve
$D$ is a smooth, nondegenerate, quartic rational curve.  
\end{thm}

\begin{proof}
This follows from the proof of ~\cite[theorem 8.1]{HRS1}.
\end{proof}

We have an induced rational transformation $g':\widetilde{H}\rightarrow
\Hdg{4}{0}(X)$ which sends a pair $(f,h)$ to the residual curve $C'$.
Also by lemma~\ref{lem-scr40}, we see that $g':\widetilde{H}\rightarrow
\Hdg{4}{0}(X)$ is (birationally) the parameter space for pairs
$(f',h')$ where $f':\Sigma\rightarrow \PP^4$ is a cubic scroll and
$h':C'\rightarrow \Sigma$ is a curve which intersects the directrix in
a degree $2$ divisor and lines of the ruling in a degree $1$ divisor,
and such that $f'(h'(C'))\subset X$ is a quartic rational curve.  Given
such a pair, the residual to $C'$ in $\Sigma\cap X$ is the quintic
elliptic $C$ we started with, so both descriptions of $\widetilde{H}$
are equivalent.

\ps

Notice that by lemma~\ref{lem-scr40}, the fiber of $g'$ over a general
point $[C']$ is the irreducible (rational) variety parametrizing pairs
$(D,\phi)$ 
where $D$ is a $2$-secant line to $C'$ and where $\phi:C'\rightarrow
D$ is an isomorphism such that $\phi$ is the identity on $C'\cap D$.  

\begin{thm}~\label{thm-irrfib2} The composite of
$g:\widetilde{H}\rightarrow \Hdg{5}{1}(X)$ and the 
Abel-Jacobi map
$u_{5,1}:\Hdg{5}{1}(X)\rightarrow J(X)$ is dominant and the general
fiber is an irreducible $6$-fold.  Therefore the Abel-Jacobi map
$u_{5,1}:\Hdg{5}{1}(X)\rightarrow J(X)$ is dominant and the general
fiber is an irreducible $5$-fold.
\end{thm}

\begin{proof}
The space of all embedded cubic scrolls $f:\Sigma \rightarrow
\PP^4$ is clearly unirational, 
in fact it is a homogeneous space for $\text{PGL}(5)$ and so it is even
rational.  So the Abel-Jacobi map is constant on the family of
complete intersections $f^*X$.  
By the residuation trick, the following two birational transformations
are pointwise (additive) inverses:
\begin{eqnarray}
\widetilde{H}\xrightarrow{g} \Hdg{5}{1}(X) \xrightarrow{u_{5,1}} J(X) \\
\widetilde{H}\xrightarrow{g'} \Hdg{4}{0}(X) \xrightarrow{u_{4,0}}
J(X). \nonumber
\end{eqnarray}
By theorem~\ref{thm-irrfib}, $u_{4,0}$ is dominant and the general
fiber is irreducible.  We have seen that $\widetilde{H}\rightarrow
\Hdg{4}{0}(X)$ is dominant and the general fiber is irreducible.  Thus
one, and hence both 
of the morphisms $\widetilde{H}_1\rightarrow J(X)$ are dominant and the
general fiber is irreducible.  Since $\widetilde{H}\rightarrow
\Hdg{5}{1}(X)$ is dominant, we conclude that $u_{5,1}$ is dominant and
the general fiber is irreducible.  
\end{proof}

\section{Double Residuation and Unirationality of the Fibers}

In the last section we introduced the space $\widetilde{H}$ which
parametrizes pairs $(f:\Sigma\rightarrow \PP^4,
h:E\rightarrow \Sigma)$ where $f:\Sigma\rightarrow \PP^4$ is a cubic
scroll and $E\subset \Sigma$ is an elliptic curve such that $f\circ
h:E\rightarrow \PP^4$ is an embedding of $E$ as a quintic elliptic
curve in $X$.  Equivalently $\widetilde{H}$ parametrizes pairs
$(f:\Sigma\rightarrow \PP^4, k:C\rightarrow \Sigma)$ where $C$ is the
rational curve residual to $E$ in $f^*X$ and $f\circ k:C\rightarrow X$
is a quartic rational curve.  In this section we will use
$\widetilde{H}$ to prove that the fibers of $u_{4,0}$ and $u_{5,1}$
are unirational.

\ps

We need a partial compactification of $\widetilde{H}$.  Let
$\overline{H}^{5,1} \subset \text{Hilb}_{5t}(X)$ denote the open
subset of the closure of
$\Hdg{5}{1}(X)$ which parametrizes Cohen-Macaulay curves (i.e. curves with
no embedded points).  Let $P(t)$ be the numerical polynomial $P(t) = 
\frac{3}{2}t^2 + \frac{5}{2}t + 1$.  Let $M~\subset~
\text{Hilb}_{P(t)}(\PP^4)$ denote the open subscheme parametrizing
connected, reduced, local complete intersection subschemes of $\PP^4$
with Hilbert polynomial $P(t)$.  The space $M$ contains an open subset
parametrizing embedded cubic scrolls, but $M$ also parametrizes mild
degenerations of embedded cubic scrolls.  We will refer to the schemes
parametrized by $M$ as \emph{generalized cubic scrolls}.  Let
$N\subset M\times \overline{H}^{5,1}$ denote the locally
closed subscheme parametrizing pairs $(\Sigma,E)$ such that $X\cap
\Sigma$ is a reduced Weil divisor and $E\subset X\cap\Sigma$.   
Notice that the Lefschetz hyperplane
theorem shows that $X$ contains no generalized cubic scrolls.  Thus
$X\cap \Sigma$ is a Cartier divisor on $\Sigma$ which contains $E$.
By ~\cite[corollary 2.7]{HRS1} the family of
residual curves to $E\subset X\cap \Sigma$ is a flat family of
Cohen-Macaulay curves.
The residual divisor $C\subset \Sigma$ is a connected, arithmetic
genus 0 curve of degree 4 contained in $X$, i.e. a point of
$\overline{H}^{4,0}$.  So we have 2 projection morphisms
$p_1:N\rightarrow \overline{H}^{5,1}$ and
$p_2:N\rightarrow \overline{H}^{4,0}$.  

\ps

One thing to notice is that if $E$ is a quintic elliptic curve whose
span is all of $\PP^4$, then $E$ is not contained in any hyperplanes,
any quadric surfaces, or in any cone over a twisted cubic.  Thus for
every point $(\Sigma,E)\in N$, we have that $\Sigma$ is a
smooth cubic scroll.  

\ps

Now define $N^2$ to be the fiber product
$N\times_{\overline{H}^{5,1}}N$.  We define the two maps $i=1,2$, 
$p_{2,i}:N^2 \rightarrow \overline{H}^{4,0}$ to be the
compositions of the projection $p_i:N^2\rightarrow N$ with
$p_2:N\rightarrow \overline{H}^{4,0}$.  

\ps

\begin{lem}~\label{lem-urfib}  For a general quartic rational curve
$[C]\in \Hdg{4}{0}(X)$, 
the fiber of $p_{2,i}:N^2\rightarrow \overline{H}^{4,0}$ is
unirational of dimension 4.
\end{lem}

\begin{proof}
Let $C\subset X$ be a general quartic rational curve.  Let $T\subset
\PP^4$ be the threefold swept out by all 2-secant lines and tangent
lines to $C$.  Then $T\cap X$ is a surface.  Let $U\subset X$ denote
the complement of $T$.  We construct a rational transformation $\rho:C
\times U \rightarrow N^2$ as follows.  For each point $p\in C$ and
$q\in U$, there is a line $L=\text{span}(p,q)$.  The quotient
projective space $\PP^4/L$ is isomorphic to $\PP^2$.  Since $L\cap C =
\{p\}$, we see that the image of $C$ under the rational projection
$\PP^4 \rightarrow \PP^2$ is a singular cubic curve.  Let $V\subset
C\times U$ denote the open set such that this singular cubic curve is
a nodal cubic.  Then the node corresponds to the unique 2-secant line
$D\subset \PP^4$ to $C$ which intersects $L$.  And there is a unique
isomorphism $\pi:C\rightarrow D$ such that $\pi(r)=r$ for each $r\in
C\cap D$ and such that $\pi(p) = s$ where $D\cap L = \{s\}$.  By
subsection~\ref{subsec-scr40} the 4-tuple $(C,D,\pi,C\cap D)$
determines a cubic scroll 
$f:\Sigma\rightarrow \PP^4$ which contains $C$ and $D$.  Let $E$
denote the residual curve to $C$ in $f^*X$.  Then $(\Sigma,E)$ is a
point of $\widetilde{H}$.  Also we know that $q$ is on $E$.  The
linear system $|2q|$ on $E$ is a $g^1_2$, i.e. a degree 2 morphism
$\pi:E\rightarrow \PP^1$.  By lemma\ref{lem-equiv1}
$\pi:E~\rightarrow \PP^1$ determines
a second cubic scroll 
$\Sigma'$ containing $E$ such that $|2q|$ is just the linear system of 
intersections of $E$ with the lines of ruling of $\Sigma'$.  For
$(p,q)\in V$ we define $\rho(p,q) = ( 
(\Sigma,E), (\Sigma',E))$.

\ps

Since every $g^1_2$ on $E$ can be expressed as $|2q|$ for some $q\in
E$, it is clear that we can obtain every cubic scroll $\Sigma'$
containing $E$ simply by varying the point $q$.  Thus we conclude that
the map $\rho:V\rightarrow N^2$ dominates the fiber of
$p_{2,i}:N^2 \rightarrow \overline{H}^{4,0}$ over $[C]$.  Since
$V$ is an open subset of the product of unirational varieties $C\times
X$, we conclude that $V$ is unirational.  Since the image of a
unirational variety is unirational, we conclude that the fiber of
$p_{2,i}:N^2 \rightarrow \overline{H}^{4,0}$ over $[C]$ is unirational.
\end{proof}

Let $P\subset N^2$ denote the irreducible component whose general
member is a point $(E,\Sigma_1,\Sigma_2)$ where $E$ is a smooth quintic
elliptic and $\Sigma_1,\Sigma_2$ are cubic scrolls.
Now consider the morphism 
\begin{equation}
r:P\rightarrow
\overline{H}^{4,0}\times \overline{H}^{4,0},
r(E,\Sigma_1,\Sigma_2) = (p_2(E,\Sigma_1),p_2(E,\Sigma_2))
\end{equation}
i.e. the pair of residual quartic curves to $E$ in $\Sigma_1\cap X$
and $\Sigma_2 \cap X$.

\begin{thm}~\label{thm-1dfib}
 For a general pair $(C_1,C_2)$ of smooth quartic rational
curves in the image of $r$, the fiber is 1-dimensional.  And the
general fiber of the Abel-Jacobi map 
\begin{equation} 
u_{4,0}:\Hdg{4}{0}(X)\rightarrow J(X)
\end{equation}
it one of the irreducible, unirational 3-folds
$p_{2,2}(p_{2,1}^{-1}([C]))$ for some $[C]\in\Hdg{4}{0}(X)$.
\end{thm}

\begin{proof}

Let $\Pi\subset \PP^4$ be a hyperplane such that $Y=\Pi\cap X$ is a
smooth cubic surface.  For a general point $p\in Y$ we can
find a quartic elliptic curve $B\subset Y$ with $p\in B$: in the
model of $Y$ as the blow-up of $\PP^2$ at 6 points, $B$ is
the proper transform of a plane cubic passing through 5 of those
points and the additional point $p$.  Also for a general point
$p\in Y$, we can find a line of type I, $L\subset X$, such that
$L\cap \Pi =\{p\}$.  Define $E=B\cup L$, so $E$ is a connected, nodal
curve of arithmetic genus 1 and degree 5.  

\begin{Claim}\label{claim-a1}
The curve $E$ satisfies the conditions in
~\cite[lemma 2.3]{HRS1},
i.e. 
\begin{equation}
H^1(L,N_{L/X}(-p))=H^1(B,N_{B/X})=0.
\end{equation}  
Therefore $\HI{5t}{X}$ is
smooth at $[E]$ and deformations of $E$ smooth the node at $p$.
\end{Claim}

We have the following short exact sequences of coherent sheaves:
\begin{equation}
\begin{CD}
0 @>>> N_{B/X} @>>> N_{E/X}|_{B} @>>> \OO_p @>>> 0 \\
0 @>>> N_{L/X} @>>> N_{E/X}|_L @>>> \OO_p @>>> 0 \end{CD}
\end{equation}
We also have a short exact sequence:
\begin{equation}
\begin{CD}
0 @>>> N_{B/Y} @>>> N_{B/X} @>>> N_{Y/X}|_{B} @>>> 0 \end{CD}
\end{equation}
The self-intersection of a quartic rational curve on a smooth cubic
scroll is $4$.  And $N_{Y/X}|_{B} = \OO_{\Pi}(1)|_E$ is also degree
4.  So by Riemann-Roch we conclude that
$H^1(B,N_{B/Y})=H^1(B,N_{Y/X}|_{B}) = 0$.  Applying the long exact
sequence in cohomology to our last short exact sequence, we conclude
that $H^1(B,N_{B/X}) = 0$.  Applying the long exact sequence in
cohomology to our first short exact sequence above, we conclude that
$H^1(B,N_{E/X}|_{B}) = 0$.  

\ps

Twisting the second exact sequence above by $\OO_L(-p)$ yields an
exact sequence:
\begin{equation}
\begin{CD}
0 @>>> N_{L/X}(-p) @>>> N_{E/X}|_L(-p) @>>> \OO_p @>>> 0 \end{CD}
\end{equation}
By assumption $N_{L/X} \cong \OO_L \oplus \OO_L$, thus $N_{L/X}(-p)
\cong \OO_L(-1)\oplus \OO_L(-1)$.  In particular, $H^1(L,N_{L/X}(-p))
= 0$.  Applying the long exact sequence in cohomology to this short
exact sequence, we conclude that $H^1(L,N_{E/X}(-p)) = 0$.  Finally,
we have the short exact sequence:
\begin{equation}
\begin{CD}
0 @>>> N_{E/X}|_L(-p) @>>> N_{E/X} @>>> N_{E/X}|_{B} @>>> 0 \end{CD}
\end{equation}
Applying the long exact sequence in cohomology to this short exact
sequence, we conclude that $H^1(E,N_{E/X}) = 0$.  Thus $[E]\in
\HI{5t}{X}$ is unobstructed.  This finishes the proof of
claim~\ref{claim-a1}. 

\ps

For a general line $M\subset \Pi$ containing $p$, the residual to $L$
in $\text{span}(L,M)\cap X$ is a smooth conic.  Let $M_1,M_2$ be two
such lines.
Without loss of
generality, we may also suppose that $M_1, M_2$ are 2-secant lines to
$B$: since we are free to choose $B$ a general quartic elliptic, we may
first choose $M_1, M_2$, and then choose $B$ to pass through $p$ and one of
the other two points of $M_i\cap Y$ for each of $i=1,2$.  Let $M_i\cap B
= \{p, q_i\}$ and let 
$r_i$ denote the third point of $M_i\cap Y$.  

\ps

Define $S_i' =\text{span}(L,M_i)$ and let $S_1'', S_2''\subset \Pi$ be
smooth quadric surfaces containing $B\cup M_1$ and $B\cup M_2$
respectively.  Define 
$S_1=S_1'\cup S_1'', S_2=S_2'\cup S_2''$, thus $S_1, S_2$ are each a union
of a 2-plane and a smooth quadric surface.  
By ~\cite[lemma 6.2]{HRS1}
such a surface is a specialization of a cubic
scroll, thus $[S_1],[S_2]\in \HI{P(t)}{\PP^4}$.  

\begin{Claim}~\label{claim-a2} For
$i=1,2$, we have $H^1(S_i,I_{E/S_i}N_{S_i/\PP^4})=
H^2(S_i,I_{E/S_i}N_{S_i/\PP^4})=0$ where $I_{E/S_i}$ is the ideal
sheaf of $E\subset S_i$.  
\end{Claim}

By the deformation theory argument in the proof of ~\cite[lemma
6.2]{HRS1}, 
this vanishing result implies that the
morphism $N^2\rightarrow \overline{H}^{5,1}$ is smooth at
$(E,S_1,S_2)$. 

\ps

We have a short exact sequence:
\begin{equation}
\begin{CD}
0 @>>> N_{S_i/\PP^4}|_{S_i'}(-M_i-L) @>>> I_{E/S_i}N_{S_i/\PP^4} @>>>
N_{S_i/\PP^4}|_{S_i''}(-B) @>>> 0 \end{CD}
\end{equation}
And we have the two short exact sequences:
\begin{eqnarray}
\begin{CD}
0 @>>> N_{S_i'/\PP^4}(-M_i-L) @>>> N_{S_i/\PP^4}|_{S'}(-M_i-L) @>>>
N_{M_i/S_i''}(-p) @>>> 0 \end{CD} \\
\begin{CD}
0 @>>> N_{S_i''/\PP^4}(-B) @>>> N_{S_i/\PP^4}|_{S_i''}(-B) @>>>
N_{M_i/S'}\otimes N_{M_i/S_i''}(-p-q) @>>> 0 \end{CD}
\end{eqnarray}
From the proof of lemma~\cite[lemma 6.3]{HRS1}, we know that 
\begin{equation}
N_{S_i'/\PP^4}
\cong \OO_{\PP^2}(1) \oplus \OO_{\PP^2}(1), N_{S_i''/\PP^4} \cong
\OO_{\PP^1\times \PP^1}(1,1) \oplus \OO_{\PP^1\times\PP^1}(2,2).
\end{equation}
Thus we have $N_{S'/\PP^4}(-M_i-L)\cong \OO_{\PP^2}(-1)\oplus
\OO_{\PP^2}(-1)$ and $N_{S_i''/\PP^4}(-B)\cong
\OO_{\PP^1\times\PP^1}(-1,-1) \oplus \OO_{\PP^1\times\PP^1}$.
Moreover $N_{M_i/S'}\cong \OO_{\PP^1}(1)$ and $N_{M_i/S_i''}\cong
\OO_{\PP^1}$.  Thus our two exact sequences are:
\begin{equation}
\begin{CD}
0 @>>> \OO_{\PP^2}(-1) \oplus \OO_{\PP^2}(-1) @>>>
N_{S_i/\PP^4}|_{S_i'}(-M_i-L) @>>> \OO_{\PP^1}(-1) @>>> 0 \\
0 @>>> \OO_{\PP^1\times\PP^1}(-1,-1) \oplus \OO_{\PP^1\times\PP^1}
@>>> N_{S_i/\PP^4}|_{S_i''}(-B) @>>> \OO_{\PP^1}(-1) @>>> 0 \end{CD}
\end{equation}
It quickly follows from the long exact sequence in cohomology that
\begin{equation}
H^{j>0}(S',N_{S_i/\PP^4}|_{S'}(-M_i-L)) =
H^{j>0}(S_i'',N_{S_i/\PP^4}|_{S_i''}(-B)) = 0
\end{equation}
Applying the long exact sequence in cohomology to our first short
exact sequence, we conclude that $H^{j>0}(S_i,I_{E/S_i}N_{S_i/\PP^4})
= 0$.  This proves claim~\ref{claim-a2}

\ps

We conclude that $N^2\rightarrow \overline{H}^{5,1}$ is
smooth at $(E,S_1,S_2)$.  Since $[E]$ is a smooth point of
$\text{Hilb}_{5t}(X)$ which lies in $\overline{H}^{5,1}$, we conclude
that $N^2$ is smooth at $(E,S_1,S_2)$ and the irreducible
component of $N^2$ which contains $(E,S_1,S_2)$ dominates
$\overline{H}^{5,1}$, i.e. the irreducible component is $P$.  

\ps

\begin{Claim}~\label{claim-a3}
The fiber of $r:P\rightarrow
\overline{H}^{4,0}\times \overline{H}^{4,0}$ containing $(E,S_1,S_2)$
is one-dimensional at $(E,S_1,S_2)$.
\end{Claim}
 
Of course the residual to $L$ in
$S_i'\cap X$ is a smooth conic $D_i'$ which intersects $M_i$ in
$\{q_i,r_i\}$.  And 
the residual to $B$ in $S_i''\cap X$ is a smooth conic $D_i''$ which
contains $r_i$.  Let us define $D_i=D_i'\cup D_i''$.  Then
$r(E,S_1,S_2) = 
(D_1,D_2)$.  

\ps

Now we want to determine the dimension of every irreducible component
of the fiber $\Phi=r^{-1}(D_1,D_2)$ through $(E,S_1,S_2)$.  By
lemma~\cite[lemma 6.3]{HRS1}, we may restrict our attention to the open
subset of $\Phi$ parametrizing pairs $(C,R_1,R_2)$ such that each of
$R_1,R_2$ is in $T\cup U$, i.e. $R_i$ is either a cubic scroll or
the union of a 2-plane and a quadric surface along a line.  
But for any conic $G$ in a cubic scroll $\Sigma$ with directrix
$D$ and fiber $F$, we know that $G\sim D+F$.  Given a union of 2 such
conics, $G',G''$, we see that $3H-G'-G'' \sim D+4F$.  Thus the
residual to such a curve in $\Sigma\cap X$ would be a quintic curve of
arithmetic genus 0, not arithmetic genus 1.  We conclude that for
$(C,R_1,R_2)$ in $\Phi$, we must have $R_i$ is in $T$.  

\ps

For
an open neighborhood of $(E,S_1,S_2)$ in $\Phi$, we must have that
$D_i'$ lies in the irreducible component of $R_i$ which is a 2-plane.
Now suppose given such a $(C,R_1,R_2)$.  There is a unique 2-plane
$S_i'$ which contains $D_i'$.  Thus $R_i = S_i'\cup R_i''$ for some
smooth quadric 
surface $R_i''$.  And $C=L\cup A$ for some quartic elliptic $A$.  Now
$\text{span}(A) $ contains both $\text{span}(D_1'')$ and
$\text{span}(D_2'')$.  Thus $\text{span}(A)=\Pi$.  Finally, since
$S_i'\cap \Pi = M_i$, we can only have $S_i'\cap R_i'' = M_i$.  In
particular, we conclude that $M_i$ intersects $A$ in the points
$p,q_i$.  Thus $A\subset Y=\Pi\cap X$ is a quartic elliptic curve
passing through $p,q_1$ and $q_2$.  These points impose independent
conditions on the 4-dimensional linear system of quartic elliptic
curves residual to $D_i''$.  Thus there is a 1-dimensional linear
system of $A$'s.  

\ps

For each $A$ in this pencil of quartic elliptic curves, we have $A\cup
D_i''\subset Y$ lies in $|\OO_Y(2)|$.  Since $Y$ is linearly normal
and lies in no quadric surfaces, we have $|\OO_Y(2)|=|\OO_{\Pi}(2)|$.
Thus there is a unique quadric surface $R_i''$ containing $A\cup
D_i''$.  Since $\{p,q_i,r_i\}\subset M_i\cap R_i''$, we conclude that
$M_i\subset R_i''$.  Thus $R_i=S_i'\cup R_i''$ is a surface in $T$,
and $(A,R_1,R_2)$ is a point in the fiber $\Phi$.  So the fiber of
$\Phi$ is one-dimensional at $(E,S_1,S_2)$, which finishes the proof
of claim~\ref{claim-a3}.  

\ps

By the above, $\Phi$ is a 1-dimensional irreducible variety in a
neighborhood of $(E,S_1,S_2)$.  It follows by upper semicontinuity of
the fiber dimension that the general fiber of $r:\PP
\rightarrow \overline{H}^{4,0}(X)\times \overline{H}^{4,0}(X)$ is at most
1-dimensional.  

\ps

For a general rational quartic $[C]\in\Hdg{4}{0}(X)$, it follows by
lemma~\ref{lem-urfib} that $p_{2,1}^{-1}([C])\subset N^2$ is a
unirational 4-fold.  So the subvariety
$Z_C:=p_{2,2}(p_{2,1}^{-1}([C])$ is an irreducible, unirational
variety.  Also observe that the subvarieties $Z_C$ sweep out
$\Hdg{4}{0}(X)$: given a curve $[D]\in\Hdg{4}{0}(X)$, if we choose any 
curve $[C]\in Z_D$, then also $[D]$ is contained in $Z_C$.  

\ps

Since the fiber dimension of $p_{2,i}:N^2\rightarrow \Hdg{4,0}(X)$ is
at most 1, we conclude the dimension of $Z_C$ is at least $4-1=3$.
Since $J(X)$ 
contains no unirational subvarieties, $Z_C$ is
contained in a fiber of the Abel-Jacobi map
$u_{4,0}:\Hdg{4,0}(X)\rightarrow J(X)$.  By theorem~\ref{thm-irrfib},
the general fiber of $u_{4,0}$ is irreducible of dimension $3$.
Combining this with the fact that a general point of
$\Hdg{4}{0}(X)$ is contained in an irreducible 3-fold $Z_C$, we
conclude that the general fiber of $u_{4,0}$ is one of the general
subvarieties $Z_C$, and vice versa.  As a consequence, observe that
the general fiber dimension of $u_{4,0}$ is 1 (and not 0).  This
completes the proof of the theorem.
\end{proof}

\ps

\begin{cor}~\label{cor-51fib}  The general fiber of the Abel-Jacobi
map
$u_{5,1}:\Hdg{5}{1}(X)\rightarrow J(X)$ is an irreducible, unirational
5-fold.
\end{cor}

\begin{proof} Recall from the proof of lemma~\ref{lem-urfib} that not
only did we show that $p_{2,i}^{-1}([C])$ is unirational, but we showed
that there is a dominant, generically-finite morphism $B\rightarrow
p_{2,i}^{-1}([C])$ such that $B$ is unirational and on $B$ we can 
 produce a
section $\sigma$ of the family of elliptic curves $E$ -- in fact we used
$|2\sigma|$ as the $g^1_2$'s on $E$ to produce the surface
$\Sigma'$.  Thus we have a distinguished line on each $\Sigma'$
corresponding to $|2\sigma|$.  And if $C'$ is the residual to
$E\subset \Sigma'\cap X$, this line intersects $C'$ in a distinguished
point.  What this shows is that there is a unirational variety $B$
dominating a general fiber $Z=u_{4,0}^{-1}(p)$ such that after we
base-change to $B$, we have a section $\tau$ of the family of quartic
rational curves $C'$.  Thus the family of quartic rational curves over
$B$ is a conic bundle with a section; therefore it is a
$\PP^1$-bundle.  

\ps

It is easy to see that given a quartic rational curve $C'$ along with
a point $p$, the locus of cubic scrolls containing $C$ is canonically
birational to $\text{Sym}^2(C)\times \CC^*$, which is canonically
birational to $\AAA^3$.  Thus the fiber product
$\widetilde{H} \times_{\Hdg{4}{0}(X)} B\rightarrow B$ is canonically
birational to $\AAA^3\times B\rightarrow B$.  In particular, since $B$
is itself unirational, we conclude that
$\widetilde{H}\times_{\Hdg{4}{0}(X)} B$ is also unirational.

\ps

Consider the
composite morphism:
\begin{equation}\begin{CD}
\widetilde{H}\times_{\Hdg{4}{0}(X)} B @>>> \widetilde{H} @>>>
\Hdg{5}{1}(X).
\end{CD}\end{equation}
Define $Y$ to be the image.  Because $\widetilde{H}\rightarrow
\Hdg{5}{1}(X)$ has fiber dimension, we conclude that 
\begin{equation}\text{dim}(Y)\geq \text{dim}
\widetilde{H}\times_{\Hdg{4}{0}(X)} B - 1 = \text{dim}(B) + 3 -1 = 3+3-1
= 5.\end{equation}
So $Y$ is a unirational variety of dimension 5.  Since $J(X)$ contains
no unirational varieties, we conclude that $Y$ is contained in a fiber
of $u_{5,1}$.  Since $\widetilde{H}\rightarrow \Hdg{5}{1}(X)$ is
dominant, we 
conclude that a general point of $\Hdg{5}{1}(X)$ is contained in one
of the varieties $Y$. 
Finally, by theorem~\ref{thm-irrfib2}, we see that the general fiber of
$u_{5,1}$ is an irreducible 5-fold.  Thus we have that $Y$ equals a
fiber of $u_{5,1}$, so the general fiber of $u_{5,1}$ is an
irreducible, unirational 5-fold.
\end{proof}

\section{Quintic Rational Curves}

In this section we will prove that the Abel-Jacobi map
$u_{5,0}:\Hdg{5}{0}(X)\rightarrow J(X)$ is dominant and the general
fiber is an irreducible, unirational $5$-fold.  

\

The construction we use to understand quintic rational curves is as
follows.  Define $I\subset \Hdg{5}{0}(X)\times \mathbb{G}(1,4)$ to be
the locally closed subvariety parametrizing pairs $([C],[L])$ where
$L$ is a $3$-secant line to $C$ which is not a $4$-secant line.  There
is a unique isomorphism $\phi:C\rightarrow L$ which is the identity on
$C\cap L$.  By lemma~\ref{lem-scr50}, there is a cubic scroll
$h:\Sigma\rightarrow\PP^4$ such that $L$ is the directrix of $\Sigma$,
$C\subset \Sigma$, and $\phi$ is the restriction of the projection
$\pi:\Sigma\rightarrow L$ to $C$.  The residual to $C$ in $\Sigma\cap
X$ is a quartic rational curve $C'$ (which does not intersect $L$ and which
intersects lines of the ruling in a degree $2$ divisor).  

\begin{lem}~\label{lem-irr50}
The scheme $\Hdg{5}{0}(X)$ is irreducible of dimension $10$.  The
morphism $I\rightarrow \Hdg{5}{0}(X)$ is birational.  For a general
pair $([C],[L])$ the residual quartic curve $C'$ is a smooth,
nondegenerate quartic rational curve.
\end{lem}

\begin{proof}
This follows from ~\cite[corollary 9.3, theorem 9.4]{HRS1}.  
\end{proof}  

What are the fibers of $I\rightarrow \Hdg{4}{0}(X)$?  By
lemma~\ref{lem-g1240}, the fiber over $[C']$ for $C'$ a general
quartic rational curve is simply an open subset of the $\PP^2$ which
parametrizes the 
collection of $g^1_2$'s on $C'$.  

\begin{thm}~\label{thm-urfib3}  The general fiber of
$u_{5,0}:\Hdg{5,0}(X)\rightarrow J(X)$ is an irreducible, unirational
$5$-fold.
\end{thm}

\begin{proof}
Since $I\rightarrow \Hdg{5}{0}(X)$ is birational, we will show that the
general fiber of $u_{5,0}:I\rightarrow J(X)$ is an irreducible,
unirational $5$-fold.  

\ps

Now $I$ surjects to an open subscheme of $\Hdg{4,0}(X)$ with
irreducible fibers isomorphic to open subsets of $\PP^2$.  By the
residuation trick, we know that $u_{5,0}$ is the pointwise inverse (up
to constant translation) of the composite map:
\begin{equation}\begin{CD}
I @>>> \Hdg{4}{0}(X) @> u_{4,0} >> J(X)
\end{CD}\end{equation}
We know the general fiber of $u_{4,0}$ is an irreducible $3$-fold.
Thus we conclude that the general fiber of $I_1\rightarrow J(X)$ is an
irreducible $5$-fold.

\ps

To see that the general fiber is unirational, we once again use the
fact that for general $p\in J(X)$ there is a morphism $B\rightarrow
u^{-1}_{4,0}(p)$ such that 
$B$ is unirational and the base-change to $B$ of the universal curve
over $\Hdg{4}{0}(X)$ admits a section.  Thus the base-change of the
universal curve is birational to $B\times \AAA^1$.  So the ``relative
symmetric product'' of the universal curve is birational to
$B\times\AAA^2$.  By lemma~\ref{lem-irr50}
$I\rightarrow \Hdg{4}{0}(X)$ is 
isomorphic to an open subset of the relative second symmetric product
of the universal curve.  Thus the fiber product
$B\times_{\Hdg{4}{0}(X)} I$ is birational to $B\times \AAA^2$.  Since
$B$ is unirational, so is $B\times \AAA^2$.  And
$B\times_{\Hdg{4}{0}(X)} I$ dominates the fiber over $p\in J(X)$.
Thus we conclude that the general fiber of $I\rightarrow J(X)$ is
unirational.  
\end{proof}

\section{Quintic Curves of Genus $2$}

By B\'ezout's theorem, $X$ cannot contain a plane curve of degree
$d>3$.  Thus the next case after quintic elliptic curves is quintic
curves of genus $2$. 

\ps

In this section we will show that the Abel-Jacobi map
$u_{5,2}:\Hdg{5}{2}(X)\rightarrow J(X)$ has image
$\text{image}(u_{5,2}) = u_{1,0}(F)$.  Moreover the general fiber is
irreducible and rational of dimension $8$.  The construction we
will use to prove this is as follows.  Suppose that $C\subset X$ is a
smooth quintic curve of genus $2$.  By Riemann-Roch, 
$h^0(C,\OO_{\PP^4}(1)|_C) = 4$, so $C$ is contained in a hyperplane
section $H\cap X$.  Similarly, $h^0(C,\OO_{\PP^4}(2)|_C)=9<10$, so $C$
is contained in a quadric surface $Q\subset H$.  The residual to $C$
in $Q\cap X$ is a line $L\subset X$.  Thus there is a morphism
$f:\widetilde{\Hdg{5}{1}}(X)\rightarrow F$ where
$\widetilde{\Hdg{5}{1}}(X)$ is the normalization of $\Hdg{5}{1}(X)$.

\begin{lem}~\label{lem-irr52}
The scheme $\Hdg{5}{2}(X)$ is irreducible of dimension $10$ and the
morphism $f:\widetilde{\Hdg{5}{2}}(X)\rightarrow F$ is dominant.
\end{lem}

\begin{proof}
This follows from the proof of ~\cite[theorem 10.1]{HRS1}.
\end{proof}

What is the fiber of $f$?  Given a line $L\subset X$, to specify the
curve $C$, it suffices to specify the hyperplane $H$ containing $L$
and the quadric surface $Q\subset H$ containing $L$.  The collection
of pairs 
\begin{equation}
\PP Q^\vee = \{([L],[H])\in F\times \PP^{4\vee}: L\subset H\}
\end{equation}
is a $\PP^2$-bundle over $F$.  And the collection of triples
\begin{equation}
I=\{([L],[H],[Q]): L\subset Q\subset H\}
\end{equation}
is a $\PP^6$-bundle over $\PP Q^\vee$.  By a dimension count we
conclude that the morphism $g:\widetilde{\Hdg{5}{2}}(X)\rightarrow I$
maps birationally to an open subset of $I$.  Thus the fibers of $f$
are irreducible, rational varieties of dimension $8$.

\begin{cor}~\label{cor-urfib52}
The Abel-Jacobi map $u_{5,2}$ has image a (translate of) an open
subset of 
$-u_{1,0}(F)$, and the general fiber of $u_{5,2}$
is an irreducible, rational $8$-fold.
\end{cor}

\begin{proof}
By the residuation trick we know $u_{5,2}$ is the pointwise inverse of
$u_{1,0}\circ g$.  And $u_{1,0}$ is an embedding.  Since $g$ is
dominant, we conclude that $\text{image}(u_{5,2})$ is an open subset
of $-u_{1,0}(X)$.  Since the fibers of $u_{5,2}$ equal the fibers of
$g$, we conclude the general fiber of $u_{5,2}$ is an irreducible,
rational variety of dimension $8$.
\end{proof}

\subsection{Irreducibility for $\Hdg{d}{0}(X)$}

We have avoided using the following result in the previous section,
but we present it here to mention one important corollary.

\ps

\begin{thm}~\label{thm-irrd0}  For each $d$ the space $\Hdg{d}{0}(X)$
is an irreducible,  reduced, local complete intersection scheme of
dimension $2d$.  Moreover the general point of $\Hdg{d}{0}(X)$ is an
unobstructed curve.   
\end{thm}

This is theorem 1 of~\cite{HRS}.  We will not discuss the proof here,
but we will prove a corollary that follows from
theorem~\ref{thm-irrd0} and our analysis of $\Hdg{4}{0}(X)$ and
$\Hdg{5}{0}(X)$.  Let $\widetilde{\Hdg{d}{0}}(X)$ denote the normalization
of $\Hdg{d}{0}(X)$.

\ps

\begin{thm}~\label{thm-irrd02}  For each $d\geq 4$, the Abel-Jacobi
morphism $\alpha_{d,0}:\widetilde{\Hdg{d}{0}}(X)\rightarrow J(X)$ is
dominant and 
the general fiber is irreducible.
\end{thm}

\begin{proof}  
Let $\mathcal{H}^d\rightarrow \overline{\Hdg{d}{0}(X)}$ denote the
normalization of the closure of $\Hdg{d}{0}(X)$.  We will actually
show that $\alpha_{d,0}:\mathcal{H}^d\rightarrow J(X)$ is dominant and the
general fiber is irreducible. 
Let us denote the Stein factorization of
$\alpha_{d,0}$ as follows:
\begin{equation}\begin{CD}
\mathcal{H}^d @> \beta_d >> Z_d @> \gamma_d >> J(X)
\end{CD}\end{equation}
We need to prove that $\gamma_d$ is an open immersion.  By
theorem~\ref{thm-irrd0}, we know that $Z_d$ is irreducible.  We will
prove by induction that there is a rational section
$\epsilon_d:J(X)\rightarrow Z_d$.  It then follows that $\gamma_d$ is
an open immersion.  We have already established this result in case
$d=4$ or $5$.  Therefore suppose that $d\geq 6$ and suppose that the
theorem has been proved for all integers less than $d$ (and greater
than $3$). 

\ps

Let $\Hdg{d-2}{0}(X)_{u}\subset \Hdg{d-2}{0}(X)$ denote the
open subscheme of $\Hdg{d-2}{0}(X)$ such that the corresponding curve
is unobstructed. 
Let $\pi:\mathcal{C}^{d-2}\rightarrow \Hdg{d-2}{0}(X)_{u}$ denote the
universal curve.
The points of $\mathcal{C}^{d-2}$ parametrize pairs $([C],p)$ where
$[C]\in\overline{\Hdg{d-2}{0}}(X)$ and $p\in C$.  
  Let $L\subset X$ be a general line.  Let
$U\subset \mathcal{C}^{d-2}$ denote the open subscheme parametrizing
pairs $([C],p)$ such that $p\in X\setminus L$.  We have a family of
$2$-planes $P\subset U\times \PP^4$ whose fiber over $([C],p)$ is
$\text{span}(L,p)$.  Let $D\subset U\times X$ denote the intersection
of $U\times X$ with $P$ inside $U\times \PP^4$.  Let $D_1\subset D$
denote the constant family $U\times L\subset U\times \PP^4$.  Let
$D_2$ denote the residual family of conics.  By
~\cite[corollary 2.7]{HRS1}, $D_2$ is flat over $U$.  Let $V\subset U$
denote the open subscheme parametrizing pairs $([C],p)$ such that the
corresponding conic $C'$ is smooth and such that $C$
intersects $C'$ transversally at $p$ and in no other
points.  Then over $V$ the union $\mathcal{C}^{d-2}\cup D_2\subset
V\times X$ is a flat family of curves of degree $d$ and arithmetic
genus $0$.  So there is an induced morphism $f:V\rightarrow
\HI{dt+1}{X}$.  

\ps

We will prove that $f(V)\subset \overline{\Hdg{d}{0}}(X)$.  In fact we
will prove that every point in $f(V)$ satisfies the conditions in
~\cite[lemma 2.3]{HRS1}, namely $H^1(C',N_{C'/X}(-p))=H^1(C,N_{C/X})=0$.
It then
follows that such a point is a smooth point of $\HI{dt+1}{X}$ and
also that the node smooths so that the point is in
$\overline{\Hdg{d}{0}}(X)$. 
Let $([C],p)$
be a point of $V$, and let $C'$ be the corresponding conic.  Let
$B=C\cup C'$ denote the union.  By construction $B$ is a connected
nodal curve.  We need to prove that $H^1(B,N_{B/X})=0$.  

\ps

We have the
exact sequence:
\begin{equation}\begin{CD}
0 @>>> N_{B/X}|_{C'}(-p) @>>> N_{B/X} @>>> N_{B/X}|_C @>>> 0.
\end{CD}\end{equation}
And we have two exact sequences:
\begin{equation}\begin{CD}
0 @>>> N_{C'/X}(-p) @>>> N_{B/X}|_{C'}(-p) @>>> \OO_p @>>> 0 \\
0 @>>> N_{C/X} @>>> N_{B/X}|_{C} @>>> \OO_p @>>> 0
\end{CD}.\end{equation}
We have seen that $N_{C'/X}$ is either $\OO_{C'}(1)\oplus \OO_{C'}(1)$
or $\OO_{C'}(2)\oplus \OO_{C'}$.  Thus $H^1(C',N_{C'/X}(-p))=0$.  By
the long exact sequence in cohomology we conclude that
$H^1(C',N_{B/X}|_{C'}(-p))=0$.  By assumption $H^1(C,N_{C/X})=0$,
therefore also $H^1(C,N_{B/X}|_C)=0$.  By the long exact sequence in
cohomology we conclude that $H^1(B,N_{B/X})=0$, i.e. $B$ is
unobstructed.  

\ps

We conclude that $f(V)$ is contained in the smooth locus of
$\overline{\Hdg{d}{0}}(X)$.  Thus we can factor $f$ as $g:V\rightarrow
\mathcal{H}^d$.  By additivity $\theta(B)=\theta(C)+\theta(C')$.  And by
residuation we have $\theta(C')=-\theta(L)$ (up to a fixed constant).
Thus we conclude that the composite map:
\begin{equation}\begin{CD}
V @> g >> \mathcal{H}^d @> \alpha_{d,0} >> J(X)
\end{CD}\end{equation}
equals the pointwise inverse (up to a constant translation) of the
composite:
\begin{equation}\begin{CD}
V @>>> \mathcal{H}^{d-2} @> \alpha_{d-2,0} >> J(X).
\end{CD}\end{equation}
Since $V\rightarrow \mathcal{H}^{d-2}$ has irreducible fibers, the Stein
factorization of this composite is just the Stein factorization of
$\alpha_{d-2,0}$ (or more precisely an open subscheme).  By the
induction assumption $\gamma_{d-2}:Z_{d-2}\rightarrow J(X)$ is an open
immersion.  By the universal property of the Stein factorization of
$V\rightarrow J(X)$, there is an induced morphism
$\epsilon_d:Z_{d-2}\rightarrow Z_d$.  This is a rational section of
$\gamma_d:Z_d\rightarrow J(X)$, which shows that $\gamma_d$ is an open
immersion.  Therefore the morphism
$\alpha_{d,0}:\mathcal{H}^d\rightarrow J(X)$ is dominant and the
general fiber is irreducible, as was to be shown.
\end{proof}

\bibliography{my}
\bibliographystyle{abbrv}
 
\end{document}